# Shuffled equi-$n$-squares

M. van de Vel,
Faculteit Exacte Wetenschappen
Vrije Universiteit Amsterdam, Amsterdam, The Netherlands
and Department of Mathematics,
University of Antwerp, Antwerp, Belgium

**Abstract.** A *formal n-square* is the set of positions in an $n$ by $n$ matrix. A *shuffle* of a formal $n$-square consists of independent rotations of each row and of each column. A key result turns out to be valid at least for $n \leq 34$ and $n = 37$: Each set of $n$ positions can be mapped with one shuffle onto a *transversal* of the columns. We consider two applications to *equi-n-squares* (i.e., $n$-matrices filled with digits $0, \ldots, n-1$ in equal amounts).

First, a shuffled equi-$n$-square can be seen as a torus with $n$ colors and two orthogonal layers of $n$ rings that can be rotated. Unlike *Rubik's cube,* each permutation of colored cells can be implemented with shuffles. A lower bound on the required amount of shuffling (approximately $n/2$ for modest $n$) obtains by a simple counting argument. An upper bound of $3(-1)^{n-1} + 6n$ is shown with the aid of the key result.

Our second application invokes column transversals and a process of *indirection* to produce *theoretically unpredictable* sequences of integers in shuffled equi-$n$-squares.

The key result has been achieved with out-of-the-box thinking: optimizing position sets, averaging, computations based on number partitions, rotating subsets of a regular $n$-gon apart, and the use of cyclotomic polynomials. A few intermediate results need computer assistance. These efforts also generated a variety of (partially) unsolved problems. We selected eight of these for a brief discussion based on the available theoretical and computer evidence.

# §1. Introduction and basic notions

An *equi-n-square* (Stein [21]) is an $n$ times $n$ matrix of *digits* $0..n-1$, each occurring $n$ times. The square is *latin* if each digit occurs in each row and each column. The number $s_n$ of equi-$n$-squares is a product of a multinomial coefficient "$n^2$ choose ($n$ times $n$)" with $n!$ (the latter counts for renaming the digits). Table 1 below compares $s_n$ with the number $l_n$ of latin squares[1] for $8 \le n \le 20$ and $n = 32$.

Equi-$n$-squares and their underlying formal squares are the major subjects of this paper. In this section, we explain our motives for this research and we present some definitions and concepts needed to formulate our methods and main results in a precise way.

| $n$ | 8 | 9 | 10 | 11 | 12 | 13 | 14 |
|---|---|---|---|---|---|---|---|
| $^2log(s_n)$ | 188.900 | 253.413 | 328.645 | 414.897 | 512.442 | 621.529 | 742.386 |
| $^2log(l_n)$ | 66.560 | 92.158 | 122.909 | 159.088 | 200.947 | 248.727 | 302.634 |

| $n$ | 15 | 16 | 17 | 18 | 19 | 20 | 32 |
|---|---|---|---|---|---|---|---|
| $^2log(s_n)$ | 875.225 | 1020.244 | 1177.626 | 1347.546 | 1530.166 | 1725.642 | 5121.444 |
| $^2log(l_n)$ | 362.864 | 392.004+ | 462.202+ | 539.215+ | 623.205+ | 714.324+ | 2410.45+ |

Table 1: Equi-$n$-Squares versus Latin Squares

**1.1. Motivation.** This project once started with an intuitive idea to produce seemingly unpredictable number sequences to serve in *stream cyphers* (Menezes et al [18]). To use equi-$n$-squares for this should not be surprising. There is a growing body of results using *latin squares* (or, *quasigroups*) for *error correcting codes* (Liu [15], Dénes and Keedwell [4, ch. 9]), *cryptocodes* (Dénes and Keedwell [5], Shcherbacov [20], Grosek and Sys [8]), *message authentication* (e.g., Meyer [19]), and *non-linear pseudo-random* (*noise*) *sequences* (e.g., Koscielny [12]). Such applications often require rather large latin squares with additional demands like (pseudo) randomness (see Jacobson and Matthews [9] for this), high non-associativity, or ease of representation. In contrast, there are no a priori requirements on $n$ or on equi-$n$-squares.

Our search for proofs has lead us to a variety of methods, results, and problems. Network flows are used to obtain partitions into latin $n$-sets, common to two equi-$n$-squares, and to produce partitions into latin column transversals (cf. § 2). Counting events and averaging are a recurring theme in § 3 where we handle optimal position sets, computations with number partitions, and rotating subsets of an $n$-gon apart, making good use of cyclotomic polynomials. This allows us to derive a *key result* on shuffling $n$-sets into column transversals and two major applications, described below in detail.

The present paper deals only with the combinatorial (and some algebraic) features of the subject. A discussion of cryptologic features is deferred to a separate paper [22].

**1.2. Coordinates and indirection.** The columns and rows of an $n$-matrix ($n \ge 2$) are numbered $0, \ldots, n-1$ in right-to-left and bottom-to-top order. A *position* (*cell*) is a pair of the form (column number, row number) and the set of all positions is a *formal n-square*. The digit assignment of an equi-$n$-square is seen as a *state* of its formal square: a partition of $n^2$ positions into $n$ equal-sized parts called *colors*. A set of positions is *latin* in a given state if its positions are all colored differently.

---

1. For values of $l_n$ ($n \le 10$) consult `http://oeis.org/A002860`. For $n = 11..15$ estimates were given by McKay and Rogoyski [16]. The exact value at 11 is due to McKay and Wanless [17]. For $n > 15$ we used the lower bound $(n!)^{2n} / n^{(n^2)}$, attributed to H.J. Ryser.

The process of *indirection* in an equi-$n$-square transforms a non-negative integer $x < n^n$ into a non-negative integer $y < n^n$ as follows. If $(x_{n-1} \ldots x_0)_n$ represents $x$ in the base $n$ then $y = (y_{n-1} \ldots y_0)_n$, where $y_k$ is the digit at position $(k, x_k)$ for $k = 0, \ldots, n-1$. The sequence of $n$ positions $(k, x_k)_{k=0}^{n-1}$ represents a *transversal* of the columns. It is the *graph of the function* $k \mapsto x_k$ ($k < n$). We call this an *H-graph* (*horizontal graph*) and the indirection process is specified as H-indirection. Similarly, V-indirection reads V-graphs, which are transversals of the rows, or, equivalently, function graphs oriented along the vertical axis. In the theory of equi-$n$-squares, a set of positions which is transversal to both the rows and the columns is better known as a *complete transversal* (Dénes and Keedwell [3, p. 26]; Stein [21]).

The term "indirection" is borrowed from an operator of the C programming language that turns a memory pointer into memory content. The suggested process treats an equi-$n$-square as if it were an array of $n$-digit numbers indexed by *all* $n$-digit numbers. To maintain and exploit this illusion, a price has to be paid as we shall see shortly.

**1.3. Intrinsic bias.** The ratio $l_n : s_n$ in table 1 illustrates that an equi-$n$-square is unlikely to have latin columns and rows. Indirection output may therefore fail *systematically* to have certain digits at certain positions. The defect can be analyzed as follows.

Let $1 \leq r \leq n$ and let $C_1, \ldots, C_r$ be distinct colors. The number of $n$-sets in $\cup_{i=1}^{r} C_i$ intersecting $C_i$ for each $i = 1, \ldots, r$ is obtained with the method of inclusion-exclusion (van Lint and Wilson [14]):

$$\sum_{k=0}^{r-1} (-1)^k \binom{r}{k} \binom{(r-k) \cdot n}{n}.$$

The number of $n$-sets with exactly $r$ colors is obtained by multiplying the previous amount with the binomial coefficient "$n$ choose $r$". The expected number of colors in an $n$-set is

$$E := \sum_{r=1}^{n} \frac{r \cdot \binom{n}{r} \cdot \sum_{k=0}^{r-1} (-1)^k \binom{r}{k} \binom{(r-k) \cdot n}{n}}{\binom{n^2}{n}}.$$

The expected number of missing colors in an $n$-set is $B := n - E$. Apparently, the state of a formal $n$-square doesn't matter. The number $B$, with its various interpretations (missing rows, missing columns, missing digits), may be dubbed the *intrinsic bias* of the formal $n$-square.

| $n$ | 2 | 3 | 4 | 5 | 6 | 7 | 8 | 9 | 10 | 11 | 12 | 13 | 14 | 15 | 16 | 17 |
|---|---|---|---|---|---|---|---|---|---|---|---|---|---|---|---|---|
| $E$ | 1.67 | 2.29 | 2.91 | 3.54 | 4.17 | 4.80 | 5.43 | 6.06 | 6.70 | 7.33 | 7.96 | 8.59 | 9.22 | 9.85 | 10.49 | 11.12 |
| $B$ | 0.33 | 0.71 | 1.09 | 1.46 | 1.83 | 2.20 | 2.57 | 2.94 | 3.30 | 3.67 | 4.04 | 4.41 | 4.78 | 5.15 | 5.51 | 5.88 |
| $B_{\mathbb{N}}$ | 0.50 | 0.89 | 1.27 | 1.64 | 2.01 | 2.38 | 2.75 | 3.12 | 3.49 | 3.86 | 4.22 | 4.59 | 4.96 | 5.33 | 5.70 | 6.07 |

| $n$ | 18 | 19 | 20 | 21 | 22 | 23 | 24 | 25 | 26 | 27 | 28 | 29 | 30 | 31 | 32 | 33 |
|---|---|---|---|---|---|---|---|---|---|---|---|---|---|---|---|---|
| $E$ | 11.75 | 12.38 | 13.01 | 13.66 | 14.28 | 14.91 | 15.54 | 16.17 | 16.80 | 17.44 | 18.07 | 18.70 | 19.33 | 19.97 | 20.60 | 21.23 |
| $B$ | 6.25 | 6.62 | 6.99 | 7.35 | 7.72 | 8.09 | 8.46 | 8.83 | 9.19 | 9.56 | 9.93 | 10.30 | 10.67 | 11.03 | 11.40 | 11.77 |
| $B_{\mathbb{N}}$ | 6.43 | 6.80 | 7.17 | 7.57 | 7.96 | 8.27 | 8.64 | 9.01 | 9.38 | 9.75 | 10.11 | 10.48 | 10.85 | 11.22 | 11.59 | 11.95 |

Table 2: Expected number $E$ of rows and bias $B$ of an equi-$n$-square versus $n$-digit bias $B_{\mathbb{N}}$.

Table 2 displays some computed values of $E$ and $B$, together with the expected number $B_{\mathbb{N}}$ of digits absent in a nonnegative number $< n^n$ (with leading zeros if needed). Note that $B < B_{\mathbb{N}}$: values in the base $n$, read in an equi-$n$-square from a random sequence of $n$ positions, tend to have slightly more different digits. Yet implementations with $n = 16$ or $32$ behave well in most of the demanding statistical tests of [6] on randomness provided the square is shuffled as prescribed below.

Stein [21, cor. 5.2] has shown that an equi-$n$-matrix must have a row or column with at least $\sqrt{n}$

distinct colors. This result is sharp, despite the much higher expectation of colors in a *generic n*-set.

**1.4. Shuffling.** Each row (resp., column) of a formal *n*-square is a copy of $\mathbf{n} := \mathbb{Z}/n\mathbb{Z}$ by interpreting numbers as column numbers (resp., row numbers) modulo *n*. We often identify the set $\mathbf{n}$ with $\{0, 1, \ldots, n-1\}$. A *rotation* by an integer amount *a* is a function of type $\mathbf{n} \to \mathbf{n}$ with $x \mapsto x + a \pmod{n}$ --suggesting a geometric interpretation of $\mathbf{n}$ as a *regular n-gon*. Our terminology extends to rows and columns of a formal *n*-square.

We avoid a systematically biased indirection output by permuting positions. A *V-move* (*H-move*) independently rotates each column (row) by some amount. A composition of such *elementary moves* amounts to a sequence alternating between V- and H-moves. A V-move followed by an H-move is called a *VH-shuffle*; in the opposite order, we have an *HV-shuffle*. We also use "shuffle" as a unit of measurement, referring to a composed move of type "HVHVH" as "2½ shuffles".

**Key Result.** *For $2 \leq n \leq 34$ and for $n = 37$, each set of n positions in a formal n-square can be mapped onto some H-graph (column transversal) with one VH-shuffle.*

The physical model of a shuffled equi-*n*-square is a *torus* with two orthogonal layers of *n* rotating rings and with (true) colors replacing digits. It may remind one of *Rubik's cube* (Joyner [10]) and the problem of recovering from a disturbed state. The first application handles a similar problem on the torus. Note that an adapted formulation of part (1) fails for Rubik's cube.

**Theorem A.** *Let $n \geq 2$.*
**(1)** *Every permutation of colored cells of an equi-n-square can be implemented with a composition of shuffles.*
**(2)** *For each equi-n-square there exists a second one which is at least $^2\log(s_n)/(2n \cdot {}^2\log(n))$ shuffles away.*
**(3)** *(n as in the Key Result.) Every two equi-n-squares are at most $6n \pm 3$ shuffles away ("-" for n even, "+" otherwise).*

For convenience, consider a *standard operation mode* repeating the following cycle of actions to produce a sequence of outputs from an equally long sequence of inputs.
**1.** Perform two HV-shuffles (shuffle input deliberately left unspecified).
**2.** Perform an H-indirection with a non-negative input integer $< n^n$.
**3.** Output the indirected integer.
Another major objective is to prove the following result.[2]

**Theorem B.** *(n as in the Key Result.) Let $l > 0$ be an integer. Given an equi-n-square and two sequences of l non-negative integers of size $< n^n$, the standard operation mode turns the first (input) sequence into the second (output) sequence with some sequence of $2 \cdot l$ shuffles.*

Theorem B is interpreted as *theoretical unpredictabilty* of the output, given the input. A realistic discussion of cryptologic features is presented in [22].

Section 2 largely prepares the way for Theorem A and proves parts 1 and 2. The Key Result is the conclusion of the entire section 3. Theorems A (part 3) and B are derived and discussed in section 4, closing with a documented list of eight problems.

---

**2.** A sophisticated version of Theorem B involves a pair of *orthogonal* equi-*n*-matrices. As each cell now holds a *unique digit pair*, some "tempering" of the output is necessary (taking an additional digit of information). This results into unpredictable 2*n*-digit output from $(n+1)$-digit input. See [22].

## §2. States, transitions and shuffles

A *transition* $P \to Q$ of equi-$n$-squares $P, Q$ is a bijective function $f$ of the underlying formal $n$-square into itself such that the $P$-color at $p$ equals the $Q$-color at $f(p)$ for each position $p$. This mimics a *physical permutation* of cells-with-content. Two transitions between the same pair of equi-$n$-squares are considered *equivalent,* suggesting a *groupoid* point of view where each pair of equi-$n$-squares has only one transition (morphism), *representable* by different bijections.

In this section we derive some results on the structure of states and on representing transitions, e.g., as composed shuffles. Our first result gives Theorem A (part 1) of the introduction.

**2.1. Shuffle Theorem.** *For each $n \geq 2$, every transition between two equi-$n$-squares can be represented with a composition of shuffles (either HV or VH). In fact, for $n$ even, every permutation of positions equals a composition of shuffles; for odd $n$, only even permutations are.*

**Proof.** To each position $(k, r)$ with column number $k$ and row number $r$ we assign a ranking number $k + nr$. A $n^2$-cycle, moving each position to a position that ranks one higher modulo $n^2$, can be achieved by an H-move, rotating each row one unit to the left, followed by a V-move, rotating the 0-th (i.e., rightmost) column up one unit, leaving all other columns at rest.

We next describe how to implement a transposition of the positions ranked 0 and 1 as a composed HV-shuffle. For $n = 2$, this can be done with an H-move, rotating the bottom row one unit to the left. For $n > 2$ odd, there is a potential problem as each shuffle is an even permutation: no sequence of shuffles can produce *exactly* a transposition. This is where equivalence comes into play. The next method applies regardless of the parity of $n$. If the positions 0 and 1 have the same color, the transposition is equivalent with the identity. If the colors are different, we can find a position ranked $x \geq n$ with the color of position 1. Then $x = k + nr$ with $r > 0$ and the 3-cycle (0 1 $x$) is equivalent with (0 1). It can be produced with 2½ HV-shuffles as follows (where $row(i)$ and $col(j)$ denote the $i$-th row and $j$-th column and all operations are modulo $n$):

$$row(0) + 1, row(r) + (1 - k); \quad col(1) - r; \quad row(0) - 1; \quad col(1) + r; \quad row(r) + (k - 1).$$

As any symmetric group $S_m$ is generated by the $m$-cycle (0 1 ... $m - 1$) and the transposition (0 1), the result for HV-shuffles follows. With due adaptations, the argument works for VH-shuffles too.

Elaborating an argument above, it can be shown that *any* 3-cycle can be performed *exactly* by a sequence of (at most 4½) shuffles. Hence shuffles of an equi-$n$-square generate at least the alternating group $A_{n^2}$. For $n$ even, there are odd generators and we obtain the full symmetric group $S_{n^2}$. □

The Shuffle Theorem suggests the question as to the amount of shuffling needed for a transition. We have the following lower bound, cf. Theorem A (part 2).

**2.2. Theorem.** *For each equi-$n$-square there exists a second one which is no less than*

$$d_n := \left\lceil \frac{{}^2 \log(s_n)}{2n \cdot {}^2 \log(n)} \right\rceil$$

*shuffles away ($s_n$ denotes the number of states and $\lceil x \rceil$ denotes upper integer approximation of $x$).*

Refer to table 1 for values of ${}^2 \log(s_n)$. One can verify that $d_n$ equals $\lceil n/2 \rceil$ for all odd $n \leq 100$ and all even $n < 30$; it equals $1 + n/2$ for $30 \leq n \leq 100$ even.

**Proof.** Given an equi-$n$-square $S$, assume we can reach every equi-$n$-square from $S$ with at most $d_n - 1$ shuffles. We reach at most $n^{2kn}$ different equi-$n$-squares with 1-1 funtions composed of $k$ shuffles. However, $n^{2kn} < s_n$ if $k = d_n - 1$, a contradiction. Hence either some square cannot be reached from $S$ at all (contradicting thm. 2.1) or some square can be reached only with $d_n$ or more shuffles. □

For a first guess of an *upper* bound we may use our proof of thm. 2.1. Count 2½ shuffles for each transposition (0 1) and one (HV-)shuffle for each long cycle (note that *powers* of the long cycle can be implemented with one shuffle, too). With these two cycles, generate all transpositions, then general cycles. This leads to an estimate of $O(n^3)$ shuffles for a generic transformation.

We can do a lot better with the results below, in combination with what is achieved in section 3. The actual upper bound theorem and its proof are postponed to section 4.

**2.3. Lemma** *Let $0 \leq l < k$ and $0 < n$ be integers. Then any two partitions of a set of size $nk - l$ into $n$ parts of size $\leq k$ have a common transversal avoiding $k - l - 1$ chosen positions.*

**Proof.** Let $U, W$ be two partitions of a set $V$ into $n$ parts of size $\leq k$. We model the situation with a network as follows. The vertex set consists of a source $s$, a sink $t$, and the three collections $U, V, W$ (formally assumed disjoint). Each part in $U$ ($W$) has an outgoing (incoming) arrow to (from) each member of $V$ which it contains. There is an arrow from $s$ to each part of $U$ and from each part of $W$ to $t$. All arrows have capacity 1. Given a set $X$ of vertices, we let $X_U, X_V, X_W$ denote the intersection of $X$ with $U, V, W$, respectively.

Obviously, there is a cut of capacity $n$ between $s$ and the remaining vertices. Let $(S, T)$ be any cut of the network with $s \in S$ and $t \in T$. Given $S' \subseteq S$ and $T' \subseteq T$, we denote by $S' \to T'$ the set of all arrows from a vertex in $S'$ to a vertex in $T'$. Expressions of type $\{s\} \to X$ and $X \to \{t\}$ are shortened to $X$ . A set name is taken to stand for its cardinality if the context requires a number. We wish to show that the capacity $S \to T$ of the cut is at least $n$:

(*)    $T_U + (S_U \to T_V) + (S_V \to T_W) + S_W \geq n$.

We may assume that $T_U + S_W < n$. Members of $S \setminus \{s\}$, resp., of $T \setminus \{t\}$, will be referred to as *used*, resp., *unused*, positions or parts. Clearly, $S_U \to T_V$ is at least the number of unused positions minus $k$ times the number of unused parts of $U$. Also, $S_V \to T_W$ is at least the number of used positions minus $k$ times the number of used parts of $W$. Adding up the two inequalities, we find that the middle terms at the left of (*) count for at least $k \cdot n - l - k \cdot (T_U + S_W)$. After adding $T_U + S_W$, we see that the left side of (*) is at least $n + k - 1 - l \geq n$.

By the theorem of Ford and Fulkerson [14] there is an integer flow of strength $n$. The flow uses a subset of $V$ of size $n$ which is a transversal of both partitions. Omitting any $k - l - 1$ chosen positions voids no part, yielding a transversal avoiding the chosen positions. □

This yields part (1) of the next result by natural induction.

**2.4. Corollary** (Structure of transitions). *Let $P, Q$ be equi-n-squares.*
**(1)**   *Given an integer $k$ with $1 \leq k \leq n$ and a set $V$ of $kn$ positions containing each color $k$ times in either state, there is a partition of $V$ into sets of size $n$ which are latin both in $P$ and in $Q$.*
**(2)**   *Given a set $W$ of positions consuming the same amount of each $P$-color and a representation $g$ of the transition $P \to Q$ with $g(W) = W$, there is a partition of the remaining positions into sets of size $n$, together with a representation $f$ of $P \to Q$ that agrees with $g$ on $W$ and maps each part onto itself.*

As to (2), note that $W$ (hence also its complement $V$) consumes equal amounts of $P$-colors and of $Q$-colors. We apply (1) on $V$, yielding parts of size $n$ which are latin in both states. Clearly, we have a representation of $P \to Q$ wich equals $g$ on $W$ and is defined on a part $S$ by assigning to a position (as colored in state $P$) the position in $S$ with the same color in state $Q$.

In latin *n*-squares, a collection of *kn* positions with each row, column, and color occurring *k* times is usually called a *k-plex*. Contrasting with the above result, there exist latin *n*-squares decomposable into *k*-plexes with $n = 2mk$, such that none of the *k*-plexes contains a *l*-plex for $1 \leq l < k$ (see [2]).

For a different application of lemma 2.3, we combine one coloring with the column (or row) partition of a formal *n*-square.

**2.5. Corollary** (Structure of equi-*n*-squares). *An equi-n-square can be partitioned into n latin H-graphs (resp., V-graphs).* □

Although by lemma 2.3 one can avoid selected positions in producing latin graphs, one cannot *force* two selected positions into one latin graph. For each $n > 2$ we found an example of an equi-*n*-square with two positions in a different row, column, and color, which do not fit together in a latin (V or H) graph.

Cor. 2.5 may suggest a strengthening of cor. 2.4(2) involving a partition into *latin graphs shared by two states*. However, with two rows colored RWB and RBW in the first state and RWR, WBB in the second state, no commonly latin graph is found in these two rows. This already shows that the strengthened prop. 2.4(1) fails for each $n \geq 3$. If the remaining row is colored RWB in both states, we obtain two equi-3-squares with no common partition into latin graphs.

It is also tempting to conclude from cor. 2.5 the structural result that there be a "network" of *n* latin H-graphs and *n* latin V-graphs, each kind partitioning the underlying formal square and with each H-graph meeting each V-graph in one position. Without the "latin" condition on one partition, such networks exist by lemma 2.3. However, if we start with a partition into latin H-graphs, the desired additional V-graphs must be latin complete transversals in a modified square, obtained from the original one by interpreting the rows, columns, and colors as, respectively, the H-graphs, the original colors, and the original rows. This relates with the *Brualdi-Ryser-Stein conjecture* (Brualdi and Ryser [1]). An example of Stein [21] shows that for any $n \geq 2$ an equi-*n*-square (even one with latin rows) need not possess a latin transversal of size *n*. Refer to question Q3 in 4.7.

The next counting result is related with the subject of section 3. All shuffles are considered to be of the same type --either HV or VH-- and all function graphs must have the same type (H or V).

**2.6. Proposition.** *The average number sh(n) of shuffles mapping a given graph to a given array of n positions, preserving the indexation, satisfies $1 < sh(n) < \sqrt{e} \approx 1.64872$ (e is the Euler number). Hence, the average number of shuffles mapping a given n-set onto a given graph lies in between n! and $n!\sqrt{e}$.*

**Proof.** The number of shuffles is $n^{2n}$; the number of *n*-arrays is $(n^2)_n$, and every shuffle maps the given function graph onto a unique array in the proper way. Hence the required average is

$$sh(n) = \frac{n^{2n}}{(n^2)_n} = \frac{n^{2n}}{n^2 \cdot (n^2 - 1) \cdot \cdot (n^2 - n + 1)} > 1.$$

As to the second inequality, we show that the natural logarithm of *sh(n)* is $< \frac{1}{2}$. Equivalently,

$$-\ln((1 - \frac{1}{n^2})(1 - \frac{2}{n^2}) \cdot \cdot (1 - \frac{n-1}{n^2})) = -\sum_{k=1}^{n-1} \ln(1 - \frac{k}{n^2}) < \frac{1}{2}.$$

Note that

$$-ln(1 - \frac{k}{n^2}) = \sum_{i=1}^{\infty} \frac{k^i}{i \cdot n^{2i}} \quad (1 \leq k < n).$$

The *i*th term $t_i$ (containing $k^i$) of the sum for $k = 1 \ldots n - 1$ of these power series satisfies

$$t_i = \frac{1}{i \cdot n^{2i}} \cdot \sum_{k=1}^{n-1} k^i.$$

The rightmost sum equals a well-known polynomial expression $P_i(n)$ in $n$ (see Zwillinger [23, Chap. 1]). We use the fact that $P_{i+1}(n) < nP_i(n)$:

$$\frac{t_{i+1}}{t_i} < \frac{nP_i(n)}{(i+1)n^{2i+2}} \cdot \frac{i \cdot n^{2i}}{P_i(n)} = \frac{i}{(i+1)n}.$$

Therefore the sum of the series is at most $t_1/(1 - 1/n) < 1/2$.

As to the last statement, different orders on an $n$-set produce disjoint sets of shuffles connecting the given graph to the ordered set (array). This leads to the desired expectation. □

We verified that $sh(n)$ increases monotonically from $4/3$ at $n=2$ to $1.64871$ at $n=50,000$. Hence the upper bound of $\sqrt{e}$ is probably sharp.

## §3. Moving a set of positions into a function graph

In this section we consider the strategic problem of transforming any set of $n$ positions of a formal $n$-square into a (function) graph with one shuffle. Prop. 2.6 suggests that this problem may be settled in the affirmative. With some efforts, this will be confirmed in moderate dimensions $n$ and with the appropriate combinations of shuffle type (VH, resp., HV) and graph type (H, resp., V).

Let $S$ be a set of positions in a formal $n$-square. An *S-row* (or, a row of $S$) is a matrix row $R$ such that $S \cap R \neq \emptyset$. We also use the term with reference to $S \cap R$, especially when referring to the *size of an S-row*. The collection of all $S$-rows is denoted *rows(S)*. An *S-free* row is a matrix row disjoint with $S$. A *row-unique position of S* is a position $p \in S$ whose matrix row $R$ satisfies $S \cap R = \{p\}$. The $S$-rows of size $> 1$ are referred to as *body rows* of $S$; taken together, the involved positions of $S$ constitute the (horizontal) *body of S*. All terminology can be adapted to columns as well; we use *cols(S)* for the collection of $S$-columns.

**3.1. Optimal and weakly optimal sets.** A set $S$ of positions is *V-optimal* (*H-optimal*) if no V-move (H-move) increases the number of $S$-rows ($S$-columns). We shall concentrate on V-optimality, the results being similar for H-optimality. A set $S$ is *weakly V-optimal* if the number of $S$-rows cannot be raised by rotating *one* column of $S$. The property is strictly weaker than optimality if $n \geq 7$. To V-optimize a set of size $n$ with $c \leq n/2$ body columns (numbered $0..c-1$, say), we basically have to rotate all body columns except (say) the zero-th one. We consider each integer from 1 up to $n^{c-1}$ (not included), rotate the $i$-th column with the $i$-th digit ($i = 1, \ldots, c-1$), and check for improving the rows score. The column-unique positions of $S$ are rotated afterwards into different free rows. For $n \geq 7$, our optimization has a worst case complexity of $O(n^c)$ with $c$ up to $\lfloor n/2 \rfloor$. Refer to problem Q4 in § 4.7.

In the formulation and proof of the next proposition, a rotation $v$ of a column $K$ is identified with the $V$-move that just rotates $K$ by $v$. The number of elements in a set $A$ is denoted by $\#A$.

**Proposition 3.2.** *Let S be a weakly V-optimal set of positions in a formal n-square and let K be an S-column. Consider $S_0 = S \cap K$ with its subset $U_0$ of row-unique positions of S in K, the set F of positions in K on S-free rows, and the set $F_0$ of all remaining positions in K. Then, with $s_0 = \#S_0$, $u_0 = \#U_0$, $f_0 = \#F_0$, and $f = \#F$, we have*

$$s_0 \cdot f \leq (n - s_0) \cdot u_0; \quad s_0 \cdot f_0 \geq (n - s_0) \cdot (s_0 - u_0).$$

*In either inequality, the non-negative difference between both sides equals the sum, taken over all rotations v of K, of the deficits $\#rows(S) - \#rows(v(S))$.*

**Proof.** As to the first inequality, the column $K$ is divided among $F, U_0, S_0 \setminus U_0, F_0$. Note that the last two sets involve *exactly* the rows occupied by $S \setminus K$. Hence, if $v$ is a rotation of $K$, then the rows of $v(S)$ hit $K$ exactly in a union of two disjoint sets

$$v(S_0) \cap (F \cup U_0) \quad \text{and} \quad F_0 \cup S_0 \setminus U_0,$$

where the second is independent of $v$. We regard $K$ as a copy of **n** by row number. The pair $(p, q) \in S_0 \times (F \cup U_0)$ is seen as an *event caused by* a rotation $v$ of $K$ provided $v(p) = q$. Grouping events by their (unique) cause divides the set $S_0 \times (F \cup U_0)$ into $n$ classes. The number of events caused by a rotation $v$ is $\#(v(S_0) \cap (F \cup U_0))$. The identity produces $u_0$ events. No rotation of $K$ can produce more, as this would raise the number of $S$-rows. Hence

$$s_0 \cdot (f + u_0) = \#(S_0 \times (F \cup U_0)) = \sum_v \#(v(S_0) \cap (F \cup U_0)) \leq n \cdot u_0.$$

The difference between the two sides of the inequality is the sum of the (non-negative) quantities

$$u_0 - \#(v(S_0) \cap (F \cup U_0)) = \#rows(S) - \#rows(v(S)),$$

taken over all rotations $v$ of $K$.

The first conclusion follows at once. The second inequality follows from the first upon noticing that the sum of the left-hand sides equals the sum of the right-hand sides. □

**Corollary 3.3.** *Let $S$ be a weakly V-optimal set of positions with $f > 0$ free rows in a formal n-square.*
**(a)** *The number $u$ of row-unique positions in an S-column of size $s$ satisfies $s - f \leq u$ and*

$$0 < \frac{s \cdot f}{n - s} \leq u.$$

**(b)** *If $\#S = n$ then $S$ has at least $n/2 + 1$ rows.*

**Proof.** Part (a) is obvious from prop. 3.2; just note that $s < n$ as $f > 0$. Part (b) follows from (a) after replacing the denominator by $n - 1$ and summing the inequality over all $S$-columns. This gives $f < u_t$ with $u_t$ the total number of row-unique positions of $S$. The number $r$ of rows satisfies $f = n - r$ and $u_t < r$ as $f \neq 0$. The result follows. □

Part (b) is obviously the best possible for $n \leq 3$. The next estimate of the number of rows is based on a different method and, in general, matches the real situation much closer.

**Proposition 3.4.** *Let $S$ be a set with $c$ S-columns of size $s_0, s_1, \ldots, s_{c-1}$, respectively, in a formal n-square. Then there is a vertical move $V$ such that the number of $V(S)$-rows is at least the outcome $r_c$ of the following recursive computation.*

$$r_0 = 0; \quad r_{i+1} = r_i + \left\lceil \frac{(n - r_i) \cdot s_i}{n} \right\rceil \quad (i = 0, \ldots, c - 1).$$

**Proof.** Let $cols(S) := \{ K_i : 0 \leq i < c \}$ with $S_i := K_i \cap S$ and $\#S_i = s_i$ for each $i$. We will show by induction on $0 \leq i < c$ that there are rotations $v_j$ of $K_j$, $0 \leq j \leq i$ such that the set

$$T_{i+1} := (S_0 + v_0) \cup \ldots \cup (S_i + v_i)$$

occupies at least $r_{i+1} \leq n$ rows. We take $v_0$ as the identity and $T_1 := S_0$ with $r_1 = s_0 \leq n$. Assume $1 \leq i < c$ with rotations $v_0, \ldots, v_{i-1}$ such that the corresponding set $T_i$ occupies at least $r_i \leq n$ rows. We mark exactly $r_i$ of them. A pair $(p, R)$, consisting of a position $p \in S_i$ and an unmarked row $R$, is considered an *event caused by* the unique rotation of $K_i$ mapping $p$ into $R$. The quantity

(*) $\quad \left\lceil \dfrac{(n - r_i) \cdot s_i}{n} \right\rceil$

equals the (up-rounded) average number of events caused by one rotation. Hence there is a rotation $v_i$

with $S_i + v_i$ occupying at least (*) many unmarked rows. Note that $r_{i+1} \leq n$.

The V-move $V$, rotating $K_i$ by $v_i$ for $0 \leq i < c$, maps $S$ to the set $V(S) = T_c$ with at least $r_c$ rows. □

It appears that the recursive computation in prop. 3.4 is sensitive to the order in which the sizes are involved (e.g., (4, 3, 2, 2) gives 8 and (4, 2, 3, 2) gives 9). We therefore define the *rows value, rval(P)*, of a partition $P$ of $n$ as the maximal output of the process, taken over all permutations of $P$.

**3.5. The algorithm** `Minrows`. We want to determine a lower bound $r = r(n)$ for the number of rows, applying to all V-optimal $n$-sets. Our method is much faster than computing the minimum of $rval(P)$ among all partitions $P$ of $n$. Given that $r = \lceil n/2 \rceil + 1$ is such a lower bound (cor. 3.3(b)), we consider partitions $P$ of $n$ with the following assumption: *P consists of the column sizes of a V-optimal n-set S with exactly $r < n$ rows.* This requires $n \geq 4$ with $c > 1$ columns.

As each $S$-column holds a row-unique position of $S$ by cor. 3.3(a), the partition's maximum, $m$, satisfies $m < r$. Putting aside one maximal part $m$, we are left with at least $\lceil (n - m)/m \rceil$ parts and at most $r - m$ parts. The algorithm skips the partitions that do not support the predicted minimum of row-unique points (cor. 3.3(a)) or the assumed number of occupied rows (prop. 3.4). A surviving *critical partition* indicates a potentially sharp lower bound $r$, which is then returned. After an unsuccessful round, $r$ is incremented for a next round.

The results $r(n)$ for $8 \leq n \leq 50$ are displayed in table 3. There is no improvement of the start value, $\lceil n/2 \rceil + 1$, for $n \leq 7$.

```
r ← ⌈n/2⌉+1
repeat
  for m from r-1 down to 2         // maximum column size
    for c from ⌈(n-m)/m⌉ up to r-m // (# of columns) - 1
      for each partition P of n-m in c parts of size ≤ m
        if ( # (row-unique points for P) > r-m ) start next partition
        P ← P plus part m
        if ( rval(P) > r ) start next partition
        else return r                 // ending Minrows
  r ← r+1
```

The algorithm `Minrows`

| $n$    | 8  | 9  | 10 | 11 | 12 | 13 | 14 | 15 | 16 | 17 | 18 | 19 | 20 | 21 | 22 | 23 | 24 | 25 | 26 | 27 | 28 | 29 |
|--------|----|----|----|----|----|----|----|----|----|----|----|----|----|----|----|----|----|----|----|----|----|----|
| start  | 5  | 6  | 6  | 7  | 7  | 8  | 8  | 9  | 9  | 10 | 10 | 11 | 11 | 12 | 12 | 13 | 13 | 14 | 14 | 15 | 15 | 16 |
| $r(n)$ | 6  | 7  | 8  | 9  | 9  | 10 | 11 | 11 | 12 | 13 | 13 | 14 | 14 | 15 | 16 | 17 | 17 | 18 | 19 | 19 | 20 | 21 |

| $n$    | 30 | 31 | 32 | 33 | 34 | 35 | 36 | 37 | 38 | 39 | 40 | 41 | 42 | 43 | 44 | 45 | 46 | 47 | 48 | 49 | 50 |
|--------|----|----|----|----|----|----|----|----|----|----|----|----|----|----|----|----|----|----|----|----|----|
| start  | 16 | 17 | 17 | 18 | 18 | 19 | 19 | 20 | 20 | 21 | 21 | 22 | 22 | 23 | 23 | 24 | 24 | 25 | 25 | 26 | 26 |
| $r(n)$ | 21 | 22 | 23 | 24 | 25 | 25 | 25 | 27 | 27 | 28 | 28 | 29 | 29 | 30 | 31 | 31 | 32 | 33 | 33 | 34 | 34 |

Table 3: Lower bound $r(n)$ on the number of rows in a V-optimal $n$-set

A useful variant of this algorithm starts at the end value $r(n)$ and returns *all* critical partitions for testing purposes.[3] This output guided us to (computer-generated) optimal $n$-sets with exactly $r(n)$ rows for $n \leq 17$ (refer to question Q5 in 4.7). Note, however, that the rows value of a partition is not always

---
**3.** The variant `Minrows` algorithm, written in well-readable Mathematica code, is given in the appendix.

accurate. E.g., for $n=16$ and $P=\{2_5, 1_6\}$ (*multiset notation:* five 2s, six 1s) we find $rval(P)=15$, but prop. 3.14 below implies that all 16-sets with a column partition of type $P$ optimize to 16 rows.

**3.6. A worked example on the sharpness of $r(n)$.** We have $r(20)=14$ and, according to the (variant) `Minrows` algorithm, the only critical partition of 20 is $\{5_4\}$. Consider a V-optimal 20-set $S$ with fourteen rows and with a column partition of this critical type. By cor. 3.3, an $S$-column with five points must have at least two row-unique points. Hence we have at most $14-4*2=6$ body rows and a body of at most $20-4*2=12$ positions. Any two $S$-columns share some $S$-row. Otherwise, we could regard them as a single column of size 10 and, as $rval(10,5,5)=15$, prop. 3.4 would contradict that $S$ is optimal. Assume *exactly* six body rows (and hence exactly eight row-unique positions). This gives a total body of exactly twelve positions and every two $S$-columns involve a common (body) row. Up to isotopism, there is only one way to achieve this, as shown in the diagram of fig. 1.

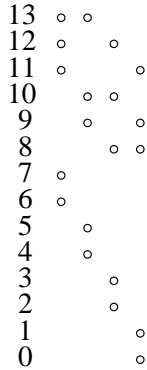

Fig. 1: The critical case in dimension 20

The inequality "$s_0 \cdot f \leq (n-s_0) \cdot u_0$" of prop. 3.2 becomes an equality with $n=20$, $f=6$, $s_0=5$, and $u_0=2$. Hence by the cited proposition, every rotation of an $S$-column results in another V-optimal configuration. By lemma 3.8 below, the leftmost $S$-column (with five members) can be rotated to avoid any three given rows, say, the ones numbered 4, 5, 13 in the diagram. This neither raises nor lowers the number of $S$-rows, so the resulting configuration is optimal again. We find that the second left $S$-column now contains three row-unique positions whereas by cor. 3.3 the other $S$-columns must have at least two. Hence we obtain an optimal configuration with *less than six body rows*.

We consider two realizations of the diagram in the formal 20-square to be *equivalent* if they differ only by permuting the order of the columns and by the location of the zero-th row. There are over $3.5 * 10^{12}$ non-equivalent realizations. We do not know whether some of these are V-optimal. Our previous argument shows that the estimate $r(20)=14$ is sharp iff there is an optimal 20-set on fourteen rows with a $\{5_4\}$ column partition and *at most five* body rows, a minor narrowing of the problem.

**3.7. The Spaghetti Effect.** The next step in achieving our key result on shuffles deals with the following problem. Given an $n$-set $S$ in a formal $n$-square with enough $S$-rows, can these rows be *rotated apart*, that is: can we find an H-move $H$ such that the rows of $H(S)$ correspond to disjoint subsets of **n**? We are speculating here on a *crumbling spaghetti effect:* a vertically stretched $n$-set has small (horizontal) sections and hence it should crumble completely under suitable "pressure" with an H-move. We need some preparatory results for a useful answer.

**Lemma 3.8.** *Let $S_1$ and $S_2$ be subsets of* **n** *with cardinality $s_1$, resp., $s_2$, and let $c = \#(S_1 \cap S_2)$. If*

(*)    $s_1 \cdot s_2 < n+c-1$,

*then $S_1$ and $S_2$ can be rotated apart in* **n**, *i.e., there is $v \in$* **n** *such that $(S_1+v) \cap S_2 = \emptyset$.*

**Proof.** Each pair $(p, q)$ with $p \in S_1$ and $q \in S_2$ is considered an *event* that is *caused by* the unique rotation moving $p$ to $q$. Note that the identity causes precizely $c$ events $(p, p)$ with $p \in S_1 \cap S_2$. Hence there are $s_1 \cdot s_2 - c$ events caused by $n-1$ non-identity rotations. By (*), one of these causes no event, and hence it maps $S_1$ outside of $S_2$. □

We next consider situations where two sets cannot be rotated apart in **n**.

**Lemma 3.9.** *Let $n \geq 2$ be integer and let $S, T$ be subsets of **n** such that $n = (\#S) \cdot (\#T)$. Then the following assertions are equivalent.*
**(1)** *$S$ and $T$ cannot be rotated apart in **n**.*
**(2)** *$S$ and $T$ do not share a positive (internal, induced) distance as subsets of the regular n-gon.*
**(3)** *$-S$ and $T$ cannot be rotated apart in **n**.*
**(4)** *For each $k \in \mathbf{n}$ there exists a unique pair $(i, j) \in S \times T$ such that $k \equiv i + j \pmod{n}$.*

**Proof.** For the equivalence (1) <--> (2) we refer to the event/cause terminology in the proof of lemma 3.8: there are $n$ events and $n$ causes. Each event $(i, j)$ (with $i \in S$ and $j \in T$) is caused by exactly one rotation. Hence some rotation causes several events iff some rotation causes none. The first statement is equivalent with the negation of (2) whereas the second is equivalent with the negation of (1).

The equivalence (2) <--> (3) follows from the previous equivalence and the fact that the mirror sets $S$ and $-S$ are isometric.

As to (3) <--> (4), assume (3) and let $k \in \mathbf{n}$. Then $(-S + k) \cap T \neq \emptyset$, which yields $i \in S$ and $j \in T$ with $-i + k \equiv j \pmod{n}$. Hence the transformation $S \times T \to \mathbf{n}$ with $(i, j) \mapsto i + j \pmod{n}$ is onto and therefore bijective. Reverse arguing yields the opposite implication. □

Let $\omega := \omega_n := e^{2\pi i/n}$ be the first primitive $n$-th root of unity. Then $\{\omega^k : 0 \leq k < n \text{ integer}\}$ gives the roots of $x^n - 1$ and represents the regular $n$-gon as a subset of the complex unit circle. The *polynomial representation* of a nonempty set $S \subseteq \mathbf{n}$ is the integer polynomial $S(x) := \sum_{i \in S} x^i$. Given $1 \leq d < n$ integer, $S$ is said to be *$d$-balanced in **n*** provided $S(\omega^d) = 0$. A set $S \subseteq \mathbf{n}$ is a *regular $m$-gon in **n*** (a regular subpolygon) if $m$ divides $n$ and $S$ is a coset of the subgroup of $m$-multiples in **n**.

The *$n$-th cyclotomic polynomial* $C_n(x)$ (Fraleigh[7, pp. 464-470]) is the minimal polynomial of $\omega_n$ over the rationals $\mathbb{Q}$ (in fact, it is in $\mathbb{Z}[x]$). In particular, $C_n(x)$ is irreducible over $\mathbb{Q}$, its roots are the primitive $n$-th roots of unity, and the product of $C_d(x)$ for $d \geq 1$ dividing $n$ equals $x^n - 1$. The degree of $C_n(x)$ is the *Euler totient* $\phi(n)$ of $n$.

Zwillinger [23, §2.3.8] has a list of $C_n(x)$ for $n \leq 30$. Packages like Maple and Mathematica have a built-in command producing cyclotomic polynomials.

**Proposition 3.10.** *Throughout, $n \geq 2$ is integer and $S, T$ are subsets of **n**.*
**(1)** *Let $n = p^2$ with $p$ prime and $\#S = p$. If $S$ is $1$-balanced in **n**, then $S$ is a regular $p$-gon within **n**.*
**(2)** *Let $n = (\#S)(\#T)$ and $1 \leq d < n$. If $S$ and $T$ cannot be rotated apart, then at least one of $S, T$ is $d$-balanced.*
**(3)** *Let $n = p \cdot q$ with $p \neq q$ prime and $\#S = p$, $\#T = q$. If $S$ and $T$ cannot be rotated apart, then one of $S, T$ is a regular subpolygon of **n** whereas the other is a complete transversal of its cosets partition. (This includes the case where both $S, T$ are regular subpolygons of **n**.)*
**(4)** *Let $n = s^2$ and $\#S = \#T = s$. If $S$ and $T$ cannot be rotated apart, then for each prime $p$ dividing $s$, one of $S, T$ is not $n/p^2$-balanced. In particular, for any three sets of size $s$ in **n**, some two can be rotated apart.*

**Proof of (1).** We may assume the set $S$ to be positioned with $0 \in S$ and the largest gap between successive points of $S$ occurring at the step reaching $0$. As the average gap between successive vertices of $S$ along the circle is exactly $p$, the degree of $S(x)$ is at most $n - p$. As $C_n(x)$ is the minimal polynomial of $\omega$ and $S(\omega) = 0$, it divides $S(x)$. Both polynomials are monic and the degree of $C_n(x)$ is $\phi(n) = n - p$, so $S(x)$ must have degree $n - p$ and hence is a regular polygon. (In fact, $S(x) = C_n(x)$ and the conclusion also follows from the formula for $C_{p^2}(x)$.)

Note: statement (1) does not extend to the case where #$S$ has two different prime factors. E.g., the 6-set $\{25, 24, 13, 12, 1, 0\}$ is 1-balanced in the 36-gon (use that $C_{36}(x) = C_{18}(x^2) = x^{12} - x^6 + 1$).

**Proof of (2).** The statement in lemma 3.9(4) can be expressed as a polynomial congruence relation,

$$S(x) \cdot T(x) \equiv \sum_{k \in \mathbf{n}} x^k \pmod{x^n - 1}.$$

Hence there is a rational polynomial $P(x)$ such that

$$S(x) \cdot T(x) - \sum_{k \in \mathbf{n}} x^k = P(x) \cdot (x^n - 1).$$

This can be rearranged as

(I) $\quad S(x) \cdot T(x) = ((x - 1) \cdot P(x) + 1) \cdot \sum_{k \in \mathbf{n}} x^k.$

For each $d$ as announced, $\omega^d$ is a root of the rightmost factor (which is $(x^n - 1)/(x - 1)$). Hence at least one of $S, T$ must be $d$-balanced.

**Proof of (3).** Part (2), applied with $d = 1, p, q$, yields that $S$ or $T$ is $d$-balanced. As $\omega, \omega^p, \omega^q$ are (respectively) a primitive $n$-th, $q$-th and $p$-th root of unity, each of the cyclotomic polynomials $C_n(x), C_q(x), C_p(x)$ divides $S(x)$ or $T(x)$. We may assume that $C_n(x)$ divides $S(x)$. Suppose that $C_q(x)$ also divides $S(x)$, say: $S(x) = C_n(x) \cdot C_q(x) \cdot Q(x)$. As $S(x)$ is monic and $C_{pq}(x), C_q(x)$ are primitive, the classical Gauss lemma (Fraleigh [7, Lemma 45.25]) yields that $Q(x)$ must have integer coefficients. Evaluating $C_n(x) \cdot C_q(x) = C_q(x^p)$ at $x = 1$ we get $q$. Taking the product with the integer $Q(1)$ must give $S(1) = p$, a contradiction. Hence $C_q(x)$ must divide $T(x)$.

If $C_p(x)$ does *not* divide $S(x)$ then $T(x)$ is divided by the product $C_p(x) C_q(x)$, which takes the value $p \cdot q$ at $x = 1$. Arguing as above, this contradicts with $T(1) = q$. So $S(x)$ is divided by $C_p(x) C_n(x)$, where the latter has degree $\phi(p) + \phi(p \cdot q) = n - q$. Arguing as in part (1), we may assume that $S(x)$ is of degree $\leq n - q$. Hence $S(x)$ must equal $C_p(x) \cdot C_n(x)$, which represents the subgroup $\{0, q, \ldots, (p-1)q\}$ of $\mathbf{n}$ (a regular $p$-subpolygon).

By lemma 3.9(4), $T$ must be a complete transversal of the coset collection of $S$ and, in fact, any such transversal will meet any rotation of $S$ (which is a coset).

Note. The general situation may be more complicated. For instance, the sets $\{25, 24, 13, 12, 1, 0\}$ and $\{10, 8, 6, 4, 2, 0\}$ share no distance and hence cannot be rotated apart in the 36-gon. Neither is a regular subpolygon.

**Proof of (4).** If $s$ happens to be a prime and if $S$ and $T$ are both 1-balanced, then both are regular $s$-gons by (1) and hence can be rotated apart in $\mathbf{n}$.

In the general situation we consider an arbitrary prime $p$ dividing $s$. Assume both sets are $n/p^2$-balanced. As $\omega^{n/p^2}$ is a primitive $p^2$-th root of unity, $S(x)$ and $T(x)$ are both divided by $C_{p^2}(x)$. Reasoning as in (2), we obtain equation (I). Its right hand side is a product that can now be divided by the square of $C_{p^2}(x)$. One factor, $\sum_{k \in \mathbf{n}} x^k$, has no multiple roots. Hence the other factor $Q(x) := (x - 1) P(x) + 1$ is divided by $C_{p^2}(x)$. However, $Q(1) = 1$ whereas $C_{p^2}(1) > 1$ since $C_{p^2}$ has only positive coefficients and at least two terms. As before, this gives a contradiction.

Given three sets $S, T, U \subseteq \mathbf{n}$ of the same size $s$, no two of which can be rotated apart, we take a prime $p$ dividing $s$ and find that exactly one of $S, T$ (say: $T$) is not $n/p^2$-balanced. Then (considering $T$ and $U$) we find by (2) that $U$ must be $n/p^2$-balanced. The pair $S, U$ contradicts the first part of (4). □

Prop. 3.10(4) will be needed near the end of this section. We now go one step further, developing conditions for rotating *multiple sets* apart. In the next proof we need the elementary fact that for any $x, y, \varepsilon \in \mathbb{R}$, if $\varepsilon > 0$ and $x + \varepsilon \leq y - \varepsilon$, then $x \cdot y < (x + \varepsilon) \cdot (y - \varepsilon)$.

**Lemma 3.11.** *Let $b \geq 2$ be an integer and for each $i = 0, \ldots, b-1$ let $S_i$ be a subset of $\mathbf{n}$ with $s_i \geq 2$ points. Let $s = s_0 + s_1 + \ldots + s_{b-1}$ and (in case $b > 2$) assume that all $s_i$ are equal or that $s - b \leq 21$. If*

(*) $\quad \dfrac{b-1}{b^2} s^2 < n,$

*then there exist rotations $r_1, \ldots, r_{b-1}$ of $\mathbf{n}$ such that $S_0 \cap \bigcap_{i=1}^{b-1} (S_i + r_i) = \emptyset$.*

**Proof.** We first consider the case $b = 2$. Using the quoted elementary fact and (*), we find that $s_0 \cdot s_1 \leq s^2/4 < n$. We may assume that the sets $S_0, S_1$ have $c > 0$ points in common. Hence $s_0 \cdot s_1 < n + c - 1$ whence by lemma 3.8, the sets can be rotated apart.

We proceed by induction on $b$. Let $b \geq 2$, assume the lemma valid for $b$ sets, and consider $S_i \subseteq \mathbf{n}$ ($i = 0, \ldots, b$) with the conditions of the lemma for $b + 1$ sets. Let $s_b$ be the minimum of the sizes $s_i$, $0 \leq i \leq b$, whence $s_b \leq s / (b+1)$. The elementary fact and the assumption (*) for $b+1$ sets yield

$$s_b \cdot (s - s_b) \leq \frac{s}{b+1} \cdot \frac{b \cdot s}{b+1} < n.$$

By lemma 3.8, *any two subsets of $\mathbf{n}$ of size $s_b$ and $s - s_b$ can be rotated apart.*

As $b + 1 > 2$, our lemma uses an additional assumption. If all $s_i$ for $0 \leq i \leq b$ are equal, then $s = (b+1) \cdot s_b$ and $s - s_b = b \cdot s_b$. We find that

$$\frac{b-1}{b^2} \cdot (s - s_b)^2 = (b-1) \cdot s_b^2 < b \cdot s_b^2 = \frac{b}{(b+1)^2} \cdot s^2 < n.$$

The alternative assumption is that $s - (b+1) \leq 21$. We have $s_b \geq 2$ and hence

(1) $\quad \dfrac{b-1}{b^2} \cdot (s - s_b)^2 \leq \dfrac{b-1}{b^2} \cdot (s-2)^2 < \left\lceil \dfrac{b}{(b+1)^2} \cdot s^2 \right\rceil \leq n,$

where the second (strict) inequality can be seen to be valid for $b \leq s - b \leq 22$.[4] Hence, with either alternative, condition (*) is available for $b$ sets $S_0, \ldots, S_{b-1}$. In case $b > 2$, we have $(s - s_b) - b \leq 21$ available to complete the requirements of the lemma for $b$ sets. Application of the induction hypothesis now yields that the sets $S_i$ ($i = 0, \ldots, b-1$) can be rotated apart. The resulting union (of size $s - s_b$) and the set $S_b$ can be rotated apart by an earlier argument. □

In more practical terms, we have the following (main) result.

**Corollary 3.12.** *Let $S$ be an n-set in a formal n-square with $b$ body rows and $f$ free rows. If $b > 2$ we assume either $f \leq 21$ or all body rows have equal size. If*

(†) $\quad \dfrac{b-1}{b^2} (b+f)^2 < n$

*(where $b + f$ equals the total body size and $b \leq f$), then $S$ can be H-moved into an H-graph.*

**Proof.** On each $S$-row, we mark the leftmost position of $S$. As the total number of $S$-rows is $n - f$, we have precizely $f$ unmarked positions in $S$, which must be located on the body rows. Hence the total number of positions on body rows is $b + f$ and, in order to have genuine body rows, it is necessary that $b \leq f$. If $b \leq 1$, a well-chosen H-move rotates the row-isolated points into different free columns. Assume $b \geq 2$. By condition (†) and lemma 3.11, we can H-move all body rows of $S$ apart. The remaining $S$-rows consist of row-unique positions, which can be rotated to different free columns. □

---

4.   For $b = 3$ the inequality (1) first fails at $s = 26$ with least $n = 127$: $(3/16) * 26^2 = 126.75 < 127$ and $(2/9) * (26-2)^2 = 128$. We expect a counterexample with four sets of sizes $2, 8, 8, 8$ where the three sets of size 8 cannot be rotated apart. Refer to the comments preceding problem Q7.

| b\f | 2 | 3 | 4 | 5 | 6 | 7 | 8 | 9 | 10 | 11 | 12 | 13 | 14 | 15 | 16 | 17 | 18 | 19 | 20 | 21 |
|---|---|---|---|---|---|---|---|---|---|---|---|---|---|---|---|---|---|---|---|---|
| 2 | 5 | 7 | 10 | 13 | 17 | 21 | 26 | 31 | 37 | 43 | 50 | 57 | 65 | 73 | 82 | 91 | 101 | 111 | 122 | 133 |
| 3 | . | 9 | 11 | 15 | 19 | 23 | 27 | 33 | 38 | 44 | 51 | 57 | 65 | 73 | 81 | 89 | 99 | 108 | 118 | 129 |
| 4 | . | . | 13 | 16 | 19 | 23 | 28 | 32 | 37 | 43 | 49 | 55 | 61 | 68 | 76 | 83 | 91 | 100 | 109 | 118 |
| 5 | . | . | . | 17 | 20 | 24 | 28 | 32 | 37 | 41 | 47 | 52 | 58 | 65 | 71 | 78 | 85 | 93 | 101 | 109 |
| 6 | . | . | . | . | 21 | 24 | 28 | 32 | 36 | 41 | 46 | 51 | 56 | 62 | 68 | 74 | 81 | 87 | 94 | 102 |
| 7 | . | . | . | . | . | 25 | 28 | 32 | 36 | 40 | 45 | 49 | 55 | 60 | 65 | 71 | 77 | 83 | 90 | 97 |
| 8 | . | . | . | . | . | . | 29 | 32 | 36 | 40 | 44 | 49 | 53 | 58 | 64 | 69 | 74 | 80 | 86 | 92 |
| 9 | . | . | . | . | . | . | . | 33 | 36 | 40 | 44 | 48 | 53 | 57 | 62 | 67 | 73 | 78 | 84 | 89 |
| 10 | . | . | . | . | . | . | . | . | 37 | 40 | 44 | 48 | 52 | 57 | 61 | 66 | 71 | 76 | 82 | 87 |
| 11 | . | . | . | . | . | . | . | . | . | 41 | 44 | 48 | 52 | 56 | 61 | 65 | 70 | 75 | 80 | 85 |
| 12 | . | . | . | . | . | . | . | . | . | . | 45 | 48 | 52 | 56 | 60 | 65 | 69 | 74 | 79 | 84 |
| 13 | . | . | . | . | . | . | . | . | . | . | . | 49 | 52 | 56 | 60 | 64 | 69 | 73 | 78 | 83 |
| 14 | . | . | . | . | . | . | . | . | . | . | . | . | 53 | 56 | 60 | 64 | 68 | 73 | 77 | 82 |
| 15 | . | . | . | . | . | . | . | . | . | . | . | . | . | 57 | 60 | 64 | 68 | 72 | 77 | 81 |
| 16 | . | . | . | . | . | . | . | . | . | . | . | . | . | . | 61 | 64 | 68 | 72 | 76 | 81 |
| 17 | . | . | . | . | . | . | . | . | . | . | . | . | . | . | . | 65 | 68 | 72 | 76 | 80 |
| 18 | . | . | . | . | . | . | . | . | . | . | . | . | . | . | . | . | 69 | 72 | 76 | 80 |
| 19 | . | . | . | . | . | . | . | . | . | . | . | . | . | . | . | . | . | 73 | 76 | 80 |
| 20 | . | . | . | . | . | . | . | . | . | . | . | . | . | . | . | . | . | . | 77 | 80 |
| 21 | . | . | . | . | . | . | . | . | . | . | . | . | . | . | . | . | . | . | . | 81 |

Table 4: Least matrix size $n(b, f)$ to rotate $b$ body rows apart, given $f$ free rows

**3.13. The Spaghetti Boundary.** Let $n(b, f)$ denote the *strict* upper integer approximation of the left side of (†). Table 4 displays the values of $n(b, f)$ for $b \leq f \leq 21$. It also shows that the restriction $f \leq 21$ of cor. 3.12 is always satisfied in dimensions $n \leq 80$.

Consider a generic set $S$ with $n$ positions in a formal $n$-square. It has $r$ rows, of which $b$ are body rows, and there are $f = n - r$ free rows. In the table, we look for values $f$ with $n(b, f) \leq n$ regardless of $b$. Sets $S$ with $n - f$ S-rows can now be H-moved into an H-graph. Given the largest such $f$ (least $r$), we find that $n - f - 1$ is an upper bound on the number of rows for a possible failure of the "crumbling spaghetti effect". The *least* upper bound is called the *Spaghetti boundary* for $n$ (taken 0 for $n = 2, 3$).

Thus table 4 yields a simple method to estimate Spaghetti boundaries for $n \leq 80$. Some of the resulting estimates (marked with a star in table 5) are down-corrected by one unit. This occurs when $n$ is a multiple of 4 with $8 \leq n \leq 28$: in the first $f$-column where the inequality $n(b, f) \leq n$ fails, it fails only for $b = n/4 + 1$. Prop. 3.14 below covers this situation and shows that the body rows can be rotated apart. The estimated Spaghetti boundary for $n = 32, 37, 43, 50$ can also be decreased by one, due to more complex reasons explained later. Refer to question Q7 in 4.7

| $n$ | 8 | 9 | 10 | 11 | 12 | 13 | 14 | 15 | 16 | 17 | 18 | 19 | 20 | 21 | 22 | 23 | 24 | 25 | 26 | 27 | 28 | 29 |
|---|---|---|---|---|---|---|---|---|---|---|---|---|---|---|---|---|---|---|---|---|---|---|
| $r(n)$ | 6 | 7 | 8 | 9 | 9 | 10 | 11 | 11 | 12 | 13 | 13 | 14 | 14 | 15 | 16 | 17 | 17 | 18 | 19 | 19 | 20 | 21 |
| $s(n)$ | 4* | 5 | 6 | 7 | 7* | 8 | 9 | 10 | 10* | 11 | 12 | 13 | 13* | 14 | 15 | 16 | 16* | 17 | 18 | 19 | 19* | 20 |

| $n$ | 30 | 31 | 32 | 33 | 34 | 35 | 36 | 37 | 38 | 39 | 40 | 41 | 42 | 43 | 44 | 45 | 46 | 47 | 48 | 49 | 50 |
|---|---|---|---|---|---|---|---|---|---|---|---|---|---|---|---|---|---|---|---|---|---|
| $r(n)$ | 21 | 22 | 23 | 24 | 25 | 25 | 25 | 27 | 27 | 28 | 28 | 29 | 29 | 30 | 31 | 31 | 32 | 33 | 33 | 34 | 34 |
| $s(n)$ | 21 | 22 | 23-1 | 23 | 24 | 25 | 26 | 27-1 | 27 | 28 | 29 | 30 | 31 | 32-1 | 32 | 33 | 34 | 35 | 36 | 37 | 38-1 |

| Starred values have been decreased by one |

Table 5: Minimal number $r(n)$ of rows of optimal sets and Spaghetti boundary $s(n)$

**3.14. Proposition.** *Given $n=4(b-1)\geq 8$, any $b$ subsets of size $2$ can be rotated apart in a regular $n$-gon.*

**Proof.** Assume $b$ sets of size 2 which cannot be rotated apart. By lemma 3.11, any $b-1$ sets among these can be rotated apart into a set $T$ of size $2(b-1)$. Let $U:=\mathbf{n}\setminus T$. If the diameter $d_0$ of the remaining 2-set would occur as a distance in $U$, we could add the corresponding 2-set to $T$ as the last rotated 2-set. Hence (counter-clock wise) rotation $r_{d_0}$ by $d_0$ maps $U$ isometrically onto $T$. In particular, the distance $d_0$ does not occur in $T$ and $r_{d_0}$ maps $T$ isometrically onto $U$ too. We deduce that the given 2-sets must all have distinct diameters. As $b\geq 3$, we may now assume $d_0\notin\{n/2,n/4\}$.

The set $T$ is partitioned into $b-1$ (rotated) original 2-sets. We pick one part $\{a_1,a_2\}$ with a diameter $d$ and we assume $r_d(a_1)=a_2$ for definiteness. There exist positions $b_i\in U$ with $r_{d_0}(b_i)=a_i$ for $i=1,2$; in particular, $r_d(b_1)=b_2$. In replacing $\{a_1,a_2\}$ by $\{b_1,b_2\}$, we obtain a set $T'$ from $T$. Like $T$, the set $T'$ is a result of rotating all given 2-sets except the one of diameter $d_0$. Hence its complement $U'$ maps isometrically onto $T'$ by $r_{d_0}$. Comparing with the previously obtained isomorphism, we see that $r_{d_0}$ maps the set $\{a_1,a_2\}$ onto $\{b_1,b_2\}$ and $r_{d_0}^2$ maps the set $\{b_1,b_2\}$ onto itself. However, $r_{d_0}^2$ is nowhere identical as $d_0\neq n/2$, nor can it swap the indices because $d_0\neq n/4$. □

We refer to question Q8 in 4.7 for additional information on rotating 2-sets apart.

Comparing the bounds in table 5 shows that $r(n)>s(n)$ for most displayed dimensions $n$, in which case the rows of a V-optimal $n$-set can be rotated apart. The corrections based on prop. 3.14 being assumed, the uncertain dimensions are $27, 30, 31, 32$ and $n\geq 35$. Combining information on the function $n(b,f)$ of Table 4 with the critical column partitions provided by the `Minrows` algorithm, we are able to settle some more cases. Our target statement on a formal $n$-square is this: *The rows of a V-optimal set with n positions in a formal n-square can be rotated apart.*

**The case $n=27$.** Our Spaghetti- and `Minrows` boundaries agree on $f=8$ free rows. According to `Minrows`, the only critical (column) partition for a V-optimal 27-set is $\{9_3\}$. This leads to at least $3*4$ row unique points by cor. 3.3(a). Hence the body size is at most $n-12=15$. According to table 4, the number $b$ of body rows can be anything between 4 and 8. In regard of the maximal body size we conclude that $b<8$. As the body is concentrated in three columns, each body row has size 2 or 3. A description of all possible body row partitions is given below (in multiset notation; cf. earlier).

| $b$ | $b+f$ | partitions of $b+f$ |
|---|---|---|
| 4 | 12 | $\{3_4\}$ |
| 5 | 13 | $\{3_3, 2_2\}$ |
| 6 | 14 | $\{3_2, 2_4\}$ |
| 7 | 15 | $\{3, 2_6\}$ |

Body row partitions for $n=27$

For $b=4$, each body row involves each of the three columns. After rotating the first three 3-sets apart, keeping the first one fixed, we obtain a 9-point set matching the fourth set at three positions. Lemma 3.8 then shows that the fourth set can also be rotated apart. Each remaining partition $(s_0, s_1, \ldots, s_{b-1})$ (with sizes $s_i$ arranged in decreasing order) can be rotated apart by using

$$(\sum_{j=0}^{i-1} s_j)\cdot s_i < n \quad (i=1,\ldots,b-1) \qquad \text{(lemma 3.8, inductively).}$$

**The case $n=30$.** A critical V-optimal set $S$ has 21 rows ($f=9$) among which are $b$ body rows, $2\leq b\leq 9$. `Minrows` throws up only one critical column partition, $\{6_5\}$, causing at least $5*3$ row-isolated points by prop. 3.3(a). Hence $S$ has at most $21-15=6$ body rows. For $b=2$ there are $b+f=11$ body positions. This is impossible since each body row of $S$ can have at most 5 points (the number of columns). We conclude that $3 \leq b \leq 6$.

| b | b+f | partitions with evaluation |
|---|---|---|
| 3 | 12 | $\{5_2, 2\}^+, \{5, 4, 3\}^+, \{4_3\}^-$ |
| 4 | 13 | $\{5, 4, 2_2\}^+, \{5, 3_2, 2\}^+, \{4_2, 3, 2\}^+, \{4, 3_3\}^-$ |
| 5 | 14 | $\{5, 3, 2_3\}^+, \{4_2, 2_3\}^+, \{4, 3_2, 2_2\}^+, \{3_4, 2\}^+$ |
| 6 | 15 | $\{5, 2_5\}^+, \{4, 3, 2_4\}^+, \{3_3, 2_3\}^+$ |

Body row partitions for $n = 30$

We prepared a small table listing all body row partitions for each value of $b$. Partitions marked with "+" satisfy the above inductive inequalities. There are two negative cases, which require more information. The partition $\{4, 3_3\}$ is the easiest. Given the number of columns (5), two body rows of respective size 3, 4 share at least two columns. We can rotate three sets of sizes 4, 3, 3 apart in a regular 30-gon. Let $T$ be the resulting 10-point set. The remaining 3-set can be rotated to share two columns with the 4-point subset of $T$. Application of lemma 3.8 yields the result.

As to the triple-4 partition, note that two body rows of this size must share 3 columns. After rotating two rows apart into an 8-set $T$, we can find two distinct rotations of the third set each causing at least three events (rotating a point of the third set into a $T$-column). We see that the remaining $30 - 2$ rotations together cause all remaining events (at most $8 \cdot 4 - 6$). So at least one of these rotations causes no event.

**The case $n = 31$.** The critical case is $f = 9$. The only critical column partition for a V-optimal $n$-set is $\{7, 6_4\}$, which gives at least $3 + 4 * 3$ row-unique positions by cor. 3.3(a). Therefore, the body size is at most 16. Table 4 shows that $3 \leq b \leq 9$. We can exclude $b \geq 8$ as this would give a body size $b + f \geq 17$. In addition, each body row lives in five columns and hence has size $\leq 5$. Here are all possible body row partitions.

| b | b+f | partitions |
|---|---|---|
| 3 | 12 | $\{5_2, 2\}, \{5, 4, 3\}, \{4_3\};$ |
| 4 | 13 | $\{5, 4, 2_2\}, \{5, 3_2, 2\}, \{4_2, 3, 2\}, \{4, 3_3\};$ |
| 5 | 14 | $\{5, 3, 2_3\}, \{4_2, 2_3\}, \{4, 3_2, 2_2\}, \{3_4, 2\};$ |
| 6 | 15 | $\{5, 2_5\}, \{4, 3, 2_4\}, \{3_3, 2_3\};$ |
| 7 | 16 | $\{4, 2_6\}, \{3_2, 2_5\}.$ |

Body row partitions for $n = 31$

In each case, inductive application of lemma 3.8 shows that the body rows can be rotated apart, except perhaps for the partition $\{4_3\}$. In this particular case, every two body rows must share at least three columns. After rotating two 4-sets apart, lemma 3.8 applies since $8 \cdot 4 < 31 + 3 - 1$.

For the next two cases $n = 32, 37$, we will show that our estimated Spaghetti boundary can be down-corrected one unit. In both cases, our (original) target statement is achieved by this.

**The case $n = 32$.** The $n(b, f)$-table 4 indicates the critical case $f = 9$ with two critical subcases: $b = 9$ and $b = 3$. If $b = 9$, there are $b + f = 18$ body positions with a $9 * 2$ body row partition. By prop. 3.14 we can rotate all body rows apart. If $b = 3$, the body has $b + f = 12$ positions and the possible body row partitions are:

$$\{8, 2_2\}, \ \{7, 3, 2\}, \ \{6, 4, 2\}, \ \{6, 3_2\}, \ \{5_2, 2\}, \ \{5, 4, 3\}, \ \{4_3\}.$$

All but the last one can be solved by inductive application of lemma 3.8. The remaining triple-4 partition is taken care of in thm. 3.15 below. We conclude that the Spaghetti boundary at $n = 32$ is at most 22, rather than 23.

**The case** $n=37$. Table 4 points at one critical case: $f=10$ with $b=3$. This gives a body size of 13 positions. Assume three body rows of sizes $s_0 \geq s_1 \geq s_2$. Then $s_0 + s_1 \leq 11$ and $s_0 \cdot s_1 \leq 30$. If $s_2 \leq 3$ then $(s_0 + s_1) \cdot s_2 = (13 - s_2) \cdot s_2 \leq 30$. If $s_2 > 3$, the partition should be $\{5, 4, 4\}$ with $(s_0 + s_1) \cdot s_2 = 36$. This allows to apply lemma 3.8 inductively. It follows that the Spaghetti boundary at $n=37$ must be estimated 26, rather than 27.

Dimensions 32, 37, 43, 50 are peculiar in table 4 because of a small peak of $n(b, f)$-values at $b=3$ in the columns $f=9,\ldots,12$. The argument for $n=37$ can be imitated to show that the Spaghetti boundary at $n=43$ is 31 rather than 32 (critical case: $f=11$, $b=3$). For $n=50$ (critical case: $f=12$, $b=3$), another similar argument shows that a decrement of the estimated Spaghetti boundary (38) depends entirely on whether three 5-sets can be rotated apart in a 50-gon. As with $n=32$, the answer is affirmative by theorem 3.15 below.

**Undecided cases** $n=35, 36, 38, 39$. For $n=35, 38, 39$, `Minrows` produces significant amounts of critical column partitions, only part of which we could handle. Dimension $n=36$ is the first with the estimated Spaghetti boundary (26) exceeding the `Minrows` boundary (25), requiring a new approach. This situation becomes permanent for $n \geq 40$.

Two questions remain: can three 4-sets (case $n=32$) or three 5-sets (case $n=50$) be rotated apart in a regular $n$-gon? The next result provides a general and affirmative answer.

**Theorem 3.15.** *If $n = t \cdot m^2$ ($m, t \geq 2$), then any $t+1$ sets of size $m$ can be rotated apart in a regular $n$-gon.*

**Proof.** Let $R_i$ for $i = 0 \ldots t$ be $m$-sets in **n** which cannot be rotated apart. Inductive application of lemma 3.8 allows us to assume that $R_i$ for $i = 1 \ldots t$ are mutually disjoint. By lemma 3.9, parts (3) and (4) we find that $(\cup_{i=1}^t R_i) - R_0 = \mathbf{n}$ (differences taken modulo $n$). Hence, as $R_i - R_0$ has at most $m^2$ elements for each $i \geq 1$, we see that

(1)   The sets $R_i - R_0$ for $i \geq 1$ partition **n** into $t$ sets of size $m^2$.

In particular, $R_0$ meets only one of $R_i$, $i \geq 1$, say: $R_1 \cap R_0 \neq \emptyset$ whereas the sets $R_0$ and $R_i$ for $i \geq 2$ are mutually disjoint. Note that $R_0$ and $R_1$ now have exchangeable roles, whence

(2)   The sets $R_i - R_1$ for $i = 0, 2, \ldots, t$ partition **n**.

We claim that $R_1 - R_0$ *is stable under addition modulo $n$.* To this end, let $v, w \in R_1 - R_0$, so $R_1 \cap (R_0 + w) \neq \emptyset \neq (R_1 - v) \cap R_0$. For any set $R \subseteq \mathbf{n}$ we have $(R - R_0) - w = R - (R_0 + w)$ (associativity), whence by (1), the sets $R_i - (R_0 + w)$ for $i \geq 1$ partition **n**. Hence as $R_1$ meets $R_0 + w$, we find

(3)   $R_i \cap (R_0 + w) = \emptyset$ for $i \geq 2$.

For any set $R \subseteq \mathbf{n}$ we have $(R - R_1) + v = R - (R_1 - v)$, whence by (2), the sets $R_i - (R_1 - v)$ for $i = 0, 2, \ldots, t$ partition **n**. Hence as $R_0$ meets $R_1 - v$, we find that

(4)   $R_i \cap (R_1 - v) = \emptyset$ for $i \geq 2$.

Suppose that $(R_0 + v + w) \cap R_1 = \emptyset$. Then $(R_0 + w) \cap (R_1 - v) = \emptyset$ and we find by (3) and (4) that the sets $R_0 + w, R_1 - v, R_2, \ldots, R_t$ are mutually disjoint. This contradicts our initial assumption and proves our claim.

As $R_1 - R_0$ is finite and contains 0, it is a *subgroup* of **n** of index $t$, which necessarily consists of the $t$-multiples among $0 \ldots n-1$ (Fraleigh [7, I§ 6]). It follows that, in particular, the internal distances of $R_0$ are $t$-folds.

The entire proof so far can be redone with any $R_i$, $i > 0$, in the role of $R_0$. Hence each of our $m$-sets exclusively has internal distances that are $t$-folds. We may therefore assume that the original sets $R_i$ are all located inside the subset $V_0$ of vertices numbered $0, t, 2t, \ldots, (m^2 - 1)t$. Let $V_i$ be the coset $V_0 + i$ ($i = 1, \ldots, t-1$). Think of copies $R'_i$ of $R_i$ for $i = 0, \ldots, t$ inside the regular $n/t$-gon, obtained by omitting the vertices of the *n*-gon which are not a *t*-fold. If, say, $(R'_0 + r) \cap R'_1 = \emptyset$ for some rotation *r* of the *n/t*-gon, then $(R_0 + t \cdot r) \cap R_1 = \emptyset$ in the *n*-gon, with both sets remaining within $V_0$. Then for each $i \geq 2$ the set $R_i$ can be rotated anywhere into the coset $V_{i-1}$, and we would effectively have rotated all given sets apart. We conclude that no two of the *m*-sets $R'_i, i = 0, \ldots, t$ can be rotated apart in $\mathbf{m^2 = n/t}$, which (as $t + 1 \geq 3$) contradicts prop. 3.10 (4). □

The previous theorem can be seen to extend prop. 3.14 (case $m = 2$). It also handles equality in condition (†) of cor. 3.12 in case $b = t + 1$, $b + f = (t+1) \cdot m$ with equal-sized body rows. Taking $n = 128$ (with $m = 8$ and $t = 2$), we see that every three sets of size 8 can be rotated apart in a 128-gon. In a footnote on the proof of lemma 3.11 we asked for an example of three 8-sets that cannot be rotated apart in a 127-gon.

## §4. The main results

In the previous section, we investigated two bounds for the number of rows of *n*-sets in a formal *n*-square: the minimal number of rows of a V-optimized set, and the maximal number of rows of a set failing to have its (body) rows rotated apart. The first exceeds the second for most dimensions $n < 40$. For such dimensions, the V-optimisation of an *n*-set can be H-optimised to an H-graph. In a few other cases where the estimated boundaries are equal, the problem was solved affirmatively using more specific information. The following (key) result summarizes our achievements so far.

**4.1. Theorem.** *Given $2 \leq n \leq 34$ or $n = 37$, each set with n positions in a formal n-square can be mapped onto some H-graph (V-graph) with one VH-shuffle (HV-shuffle).* □

In a formal *n*-square, the *expected* number of rows of an *n*-set (table 2) is rather close to --and sometimes even larger than-- the estimated Spaghetti bound for moderate *n* (table 5). This suggests that for such dimensions *n* a (near) majority of *n*-sets can be H-moved into an H-graph without a preliminary V-move. Random samples in dimensions $n = 15 \ldots 20$ revealed an *excessively large* majority of *n*-sets that can be H-moved into an H-graph. In fact, on $10^5$ random sets, the number of counterexamples decreased from 270 at $n = 15$ to merely 8 at $n = 20$. We currently have no explanation for this. For small *n*, prop. 4.2 below provides some information.

**4.2. Proposition.** *For $n = 2, 3$ each n-set in a formal n-square can be H-moved onto an H-graph and V-moved onto a V-graph. For $n = 4, 5$, each n-set satisfies at least one alternative. Both alternatives may fail for a 6-set in a formal 6-square.*

**Proof.** For $n = 2, 3$ this already appears from the `Minrows` boundary $r(n) = n$ if $n = 2, 3$. In general, a *n*-set with at most one body row (column) can be H-moved (V-moved) apart. Hence for $n = 4$, at least one of the alternatives holds in case a 4-set has *r* rows or *r* columns with $r = 1, 3, 4$. A 4-set with two rows and two columns has a $\{2, 2\}$ rows partition and its parts share a distance. The same goes for the columns partition. Hence lemma 3.8 yields both alternatives.

Consider a 5-set *S* with *r* rows. If $r \in \{1, 4, 5\}$ we have at most one body row. Suppose $r = 3$ with two body rows. The row partition is $\{2, 2, 1\}$ and the two pairs can be rotated apart by lemma 3.8. If $r = 2$, *S* will occupy at least three columns. This case is solved, *mutatis mutandis,* as with three rows.

For $n = 6$, the set $\{(0,0), (1,0), (1,1), (2,1), (0,2), (2,2)\}$ can't be moved to a graph as required. □

**4.3. Corollary.** *Let $2 \leq n \leq 34$ or $n = 37$.*
**(a)** *An array with n positions in a formal n-square can be mapped onto some H-graph with an HVH-move preserving the order. This also holds with "V-graph" and "VHV-move".*
**(b)** *Given an n-array and an H-graph (V-graph), there is an HVHV (VHVH) move mapping the array onto the graph and preserving the order.*

**Proof.** By theorem 4.1, we can map the set of positions of an array $A$ onto some V-graph $G$ by an HV-move. We map $G$ onto an H-graph $H$ by an H-move defined as follows. For each $i$ with $0 \leq i < n$ we rotate the row of position $G[i]$ so that $G[i]$ moves to the $j$-th column, where $A[j]$ is the element that we mapped to $G[i]$. Each member of $G$ ends up in a different column, whence we obtain an H-graph with the correct ordering.

As to part (b), any H-graph can be mapped orderly into any other H-graph with a V-move. □

Prop. 2.6 provides a rough "statistic" variant of cor. 4.3(b). It needs (on average) only one shuffle and is valid in all dimensions $n$. A $n$-array with $c$ columns allowing such a shuffle will allow at least $n^{n-c}$ of them. Hence the expectation of $n-c$ (the intrinsic bias, cf. table 2) suggests that such shuffles do not exist for *every* $n$-array. The simplest possible example for any $n \geq 3$ involves an array $A$ with $A[0] := (1, 0)$ and $A[1] := (0, 0)$. Regardless of the other positions of $A$, there is no shuffle (either HV or VH) mapping $A$ orderly onto the bottom row.

The two shuffles required in 4.3(b) are well-spent: the first does the optimizing and the second takes care of precision. The statement can be reformulated like this: *any 1-1 function of an n-set onto an H-graph extends to a composition of two HV-shuffles*. Combining this with both thm. 4.1 and cor. 4.3(b), we get a simple proof that --with the usual restriction on $n$-- *any 1-1 function between two n-sets in a formal n-square extends to a composition of three shuffles* (two versions: HV and VH).

We now derive the two major theorems announced in section 1. Recall that the *standard mode of operation* on an equi-$n$-square is to perform two HV-shuffles before each H-indirection.

**4.4. Theorem** (cf. Theorem B). *Let $2 \leq n \leq 34$ or $n = 37$, and let $l > 0$. Given an equi-n-square, together with two sequences, both consisting of $l$ non-negative numbers of size $< n^n$, there is a sequence of $2 \cdot l$ shuffles such that the standard mode of operation on an equi-n-square turns the first (input) sequence into the second (output) sequence.*

**Proof.** An output number can be seen as an array of $n$ digits. We can find the same array in the matrix at a suitable sequence of $n$ positions. By corollary 4.3(b), two successive shuffles can map this sequence orderly onto the target graph, described by the input number. Now the output by indirection is as desired. □

Earlier we described this result as *theoretical unpredictability of the output*. The literature provides some other combinatorial methods aiming at this goal. Knuth [11, Algorithm M on p.33] uses an index-to-value process with (standard) arrays to hide the inevitable output patterns of pseudo-random number generators. Lidl and Niederreiter [13] produce multiplexed sequences by using one m-sequence to step through another one. Koscielny [12] uses one-cell indirection in a latin square to modulate a linear m-sequence with a shifted copy.

One might consider a more general type of indirection, reading out arbitrary $n$-arrays in an equi-$n$-square. According to a remark following cor. 4.3, any two $n$-arrays in an equi-$n$-square are linked with a sequence of three shuffles (usual $n$). Hence by adapting the previous argument, this method can produce any sequence of numbers with *three shuffles* before each indirection. The method also needs a *doubled input* to determine the array to read and hence may be qualified as rather wasteful. However, with this type of indirection (explicitly allowing repeating cells) the intrinsic square bias for numbers becomes the natural bias (see § 1) and the need for shuffling decreases.

**4.5. Theorem** (cf. Theorem A, part 3). *Given $2 \leq n \leq 34$ or $n = 37$, every two equi-$n$-squares are at most $6n + 3(-1)^{n-1}$ shuffles away from each other.*

**Proof.** Let $P$ and $Q$ be equi-$n$-squares. By lemma 2.5 we can fix a latin V-graph $G$ in $P$ and a latin H-graph $H$ in $Q$. There is a 1-1 function connecting each $P$-colored position of $G$ with the position of $H$ having the same color in $Q$. We may assume a transition $f : P \to Q$ extending this 1-1 function. As in the proof of cor. 4.3(b), the inverse function $H \to G$ agrees with a VH-shuffle $\sigma$ on $H$. Hence $f$ followed by $\sigma$ keeps every position of $G$ fixed. This shows that if each transition preserving the coloring of $G$ can be implemented with $k$ (or even $k + \frac{1}{2}$) HV-shuffles for some integer $k$, then our transition $f : P \to Q$ can be implemented with $k + 1$ HV-shuffles (collapsing two successive H-moves if needed) by letting the shuffle sequence be followed by the HV-shuffle $\sigma^{-1}$. We henceforth concentrate on transitions $f : P \to Q$ leaving all positions of a latin (in $P$) V-graph $G$ invariant.

By prop. 2.4(2), we may assume a partition of $G$'s complement into $n - 1$ sets of size $n$ with $f$ mapping each part onto itself. On each part $S$ we decompose the permutation into $t(S)$ disjoint cycles $C(S)_1, \ldots, C(S)_{t(S)}$ of length $\geq 2$. Note that $t(S) \leq n/2$. The *carrier* of $f$ in $S$ can be ordered using the list of cycles and the choice of a "last point" $l(S)_i$ of $C(S)_i$ ($i = 1, \ldots, t(S)$). The given permutation of $S$ equals the inverse of the cycle $L(S) := (l(S)_1 \ldots l(S)_{t(S)})$, followed by the cycle provided by the ordered carrier.

The parts $S$ are taken in pairs with a single part remaining if $n$ is even. Assume $n$ is even; for $n$ odd, we present the necessary adaptations afterwards. On a selected pair $(S_1, S_2)$ of parts, the transition is processed as a sequence of HV-shuffles as follows.

**(Step 1)** For each $i = 1, 2$, extend the $t(S_i)$-cycle $L(S_i)$ of "last points" to an equivalent cycle of length $n/2$, e.g., by repeating one of its colors. We form an $n$-sequence in which the two resulting sequences are mixed alternatingly. Some HVHV-move maps the sequence to (say) the bottom row of the matrix (cor. 4.3(b)). The row is then rotated two positions to the right and moved back to the array with the reverse of the initial move. This makes the reverse of the product of $L(S_1)$ and $L(S_2)$ equivalent with $4\frac{1}{2}$ HV-shuffles.

**(Step 2)** This applies to the sets $S_i$ for $i = 1, 2$ separately. The ordered carrier in $S_i$ is made into an $n$-array by repeating a color. It takes two HV-shuffles to map the array onto the bottom row (cor. 4.3(b)) which is then rotated one position to the left and mapped back by reversing the initial HV shuffles. This makes the carrier cycle on $S_i$ equivalent with $4\frac{1}{2}$ HV-shuffles.

The resulting sequence of shuffles, resulting from steps (1) and (2), reduces to $12\frac{1}{2}$ HV-shuffles after collapsing successive H-moves.

The set $S$ that is not in a pair is treated with the following modifications. In (Step 1), we extend the cycle $L(S)$ to an $n$-cycle by repeating the last color and perform it directly as $4\frac{1}{2}$ shuffles with the method of (Step 2). Step (2) remains the same. This accounts for a total of $8\frac{1}{2}$ HV-shuffles.

Compiling all results, we can contract the last H-move of one shuffle sequence with the first H-move of the next shuffle sequence. This amounts to one sequence of 8 HV-shuffles, $(n-2)/2$ sequences of 12 HV-shuffles, and a leftover H-move to finish with. The latter is absorbed by a final HV-shuffle as described at the beginning of this proof, accounting for $6n - 3$ shuffles.

For $n$ odd, no set is taken single. We need a preliminary step involving *all* selected pairs $(S_{2i-1}, S_{2i})$ for $i = 1, \ldots, (n-1)/2$.

**(Step 0)** Let $(a_i, b_i)$ be the pair consisting of the last point in the representations of $L(S_{2i-1})$, resp., $L(S_{2i})$. As observed earlier, the product of these transpositions is obtained by first applying the cycle $(b_1, \ldots, b_{(n-1)/2})$ in reverse, then applying the cycle $(a_1, b_1, \ldots, a_{(n-1)/2}, b_{(n-1)/2})$. Both cycles can be performed equivalently with $4\frac{1}{2}$ HV-cycles. Hence the product of pairs, required at this step, can be achieved with $8\frac{1}{2}$ HV-shuffles.

In step (1), we mix the two given cycles as we did in the even case and perform the resulting cycle with $4\frac{1}{2}$ shuffles. Step (2) is as above. The entire state transformation is now equivalent with: the result of step (0), followed by the results of step (1), followed by the results of step (2), ended with the obligatory HV-shuffle. After contracting successive H-moves, this accounts for $6n + 3$ HV-shuffles.

□

The previous result was derived originally to get an idea of how many shuffles are needed to obtain a *generic state* of a formal square (e.g., to serve as an *initial state*). For instance, 93 shuffles should do for $n=16$. The estimates in thm. 4.5 are too high for small $n$. E.g., for $n=2$ it can be seen that any two states are connected by one HV- or VH-shuffle. They are connected by 1½ shuffle of either type. For $n \leq 5$ some economizing is possible with the aid of prop. 4.2, resulting in an estimate of $5n \pm 3$. A minor inaccuracy may occur for $n=4, 5$ due to a repeated dilemma: HVH or VHV shuffle.

**4.6. The physical equi-$n$-square.** The physical object implementing an equi-$n$-square with shuffling is a *torus* with two orthogonal layers of $n$ rings, each of which can be rotated by any integer multiple of $2\pi/n$ radians. In combination, the two layers divide the torus surface into $n^2$ rectangular faces and performing a shuffle amounts to rotating rings of one type, then rotating the other type. With faces evenly painted with $n$ colors, the object becomes a variant of *Rubik's cube* (see e.g. Joyner [10]).

Theorem 4.5 contributes to the problem of *restoring a physical state* of the torus from any other state. Estimates of our method in the *equi-octal-square* (with 64 faces, it is close in size to the standard Rubik's cube) indicate that this may take up to 45 shuffles, each shuffle requiring up to 16 rings to turn. In contrast, the cube never needs more than 20 elementary moves. One explanation might be that there are implicit restrictions on moving cells of the cube, contrasting with the Shufffle Theorem 2.1. Unsharpness of our estimates is a more probable explanation (refer to Q2 below).

**4.7. Conclusion and open problems.** We achieved two distinct goals from a combinatorial study of equi-$n$-squares and their underlying formal squares. One is a principle on which to produce unpredictable sequences of non-negative integers, the other is an upper bound on the shuffle distance between two equi-$n$-squares. Both results are valid under the common restriction $2 \leq n \leq 34$ or $n=37$ and the proofs borrow largely from the same intermediate results. Our main results, as well as some of the intermediate results, have raised various questions and we wish to discuss some of these.

The somewhat irregular restriction on $n$ derives from an interaction of several combinatorial methods (as described in § 3) and is probably not sharp. For $n$ beyond the restriction, however, little information is available. Elaboration of our methods shows that $n$-sets in a formal $n$-square are mapped by some VHV-move onto a V-graph for $n<45$ (with possible exception of $n=42$). There are usually many arrays in an equi-$n$-square representing a given $n$-digit number, offering increased opportunities to shuffle one of these arrays onto a graph. This suggests that thm. 4.4 may hold well beyond its current limitations.

Note. Our remarks following thm 4.1 indicate that performing "on purpose" 1½ shuffles before each indirection can reproduce an output sequence, obtained under standard operation mode from a given input sequence, *with a very low ratio of errors* ($15 \leq n \leq 20$). This suggests that, practically speaking, our method should work even with 1½ shuffles, alternating between H- and V-indirection.

**Q1.** *Determine all n for which thm. 4.4 holds as formulated. How many shuffles are needed for other n (if any)? Is there a finite asympotic value for $n \rightarrow \infty$?*

Our upper bound for the number of shuffles implementing a generic transition seems too high. Hard evidence has been provided above for $n \leq 5$. Our proof strategy of thm. 4.5 (realizing the transition stepwise on successive parts) necessarily includes some cleaning up of undesired changes at each step. This may suggest some further economizing. Given the considerable gap with the lower bound in thm. 2.2, the following is a challenging question.

**Q2.** *Is it true that the lower bound on the required number of shuffles implementing a generic transition is closest to reality, in other words, is this a* small world phenomenon?

The following open problems relate with our combinatorial tools. First, in regard of comments following cor. 2.5, we define a *wavy-latin n-square* to be an equi-$n$-square allowing a partition into latin

H-graphs and one into latin V-graphs with each H-graph meeting each V-graph in one position (a *wavy-latin network*). It follows from cor. 2.5 that an equi-*n*-square with latin rows or with latin columns is wavy-latin. As observed in § 2, lack of latin *n*-size transversals seems a major obstacle for an equi-*n*-square to be wavy-latin.

The defined property is invariant under color renaming, permuting the order of rows and of columns (*isotopism*), and transposition. This leads to a small collection of "types". For $n=3$, we have 9 types, of which one is not wavy-latin. For $n=4$, we counted 3 160 types of which exactly 8 are not wavy-latin.

$$\begin{array}{cc} 1 & 0 \\ 1 & 0 \end{array} \qquad \begin{array}{ccc} 1 & 2 & 1 \\ 1 & 0 & 0 \\ 2 & 2 & 0 \end{array} \qquad \begin{array}{cccc} 1 & 1 & 3 & 2 \\ 1 & 1 & 3 & 3 \\ 3 & 2 & 0 & 0 \\ 2 & 2 & 0 & 0 \end{array} \qquad \text{(non-wavy-latin } n\text{-squares)}$$

A random sample of 2 000 equi-5-squares[5] has not revealed a counterexample. In addition, all equi-5-squares with each color having a horizontal or vertical alignment of at least four cells are wavy-latin. We reduced this to 164 cases, mostly verified with computer assistence.

**Q3.** *Does the ratio (wavy-latin n-squares / equi-n-squares) tend to 1 with increasing n? Is there a size $n_0$ such that all equi-n-squares are wavy-latin for $n \geq n_0$?*

(General formal *n*-squares.) Given an *n*-set $S$ of positions, the result of counting row-unique positions in an *S*-column or counting *S*-rows may provide a sufficient reason for $S$ not to be V-optimal (props. 3.3(a) and 3.4). Otherwise, we have to rely on an algorithm of complexity $O(n^{n/2})$ for $n \geq 7$. To compute a V-optimal configuration of a given *n*-set, one could follow the computation of the rows value of its column partition (which is $O(\lfloor n/2 \rfloor!)$ at worst) to know the best order in which to rotate the individual (body) columns and the apropriate amount of the rotations. However, there is an example at $n=16$ where the rows value (15) is not the optimal value (16). See § 3.5.

**Q4.** *($n \geq 7$.) Is there a fast(er) method to optimize a generic set of n positions or to decide that a generic n-set is optimal? More specifically, is this problem NP? Is it NP complete?*

Our `Minrows` lower bound $r(n)$ for the number of rows of a V-optimal *n*-set in an equi-*n*-square is proven sharp for $n \leq 17$ by computer-supplied examples. As to the sharpness of $r(20)=14$, we found a minor reduction of the problem (cf. 3.6).

**Q5.** *Find sharp lower bounds for $r(n)$, $n \geq 18$. On replacing "V-optimal" by "weakly V-optimal", are there any n for which this gives a smaller boundary?*

The key to the major results in this paper is thm. 4.1 on mapping an *n*-set onto a function graph with one shuffle. We derived that any 1-1 function of an *n*-set onto a graph extends to a composition of two shuffles and that a 1-1 function between two *n*-sets extends to a composition of three shuffles (both VH and HV).

**Q6.** *For which n can each n-set in a formal n-square be mapped onto an H- or V-graph with one shuffle? How many shuffles are needed for other n? Is there a finite asymptotic value? Can 1-1 maps between generic n-sets be extended to a composition of less than 3 shuffles?*

Table 4 provides a rough-and-ready method to estimate the Spaghetti boundary for dimensions $\leq 80$. The method rests upon cor. 3.12, valid under the restriction (for $b \geq 3$) that the number $f$ of free rows is $\leq 21$. A close inspection of the preceding lemma 3.11 shows that this restriction is a kind of common denominator for a variety of situations. In fact, it turns out that the limit on $f$ increases with the number $b$ of body rows as described in table 6 (computed with Maple). This enlarges the range of applicability of our method. N.B.: the computed limit at $b=2$ does not contribute to lemma 3.11.

---

[5].  Considering the exception rates for $n=2, 3, 4$, a sample of 2 000 equi-5-squares is probably too small to make a solid prognosis. Unfortunately, it takes us 10 minutes on average to explore one 5-square.

| $b=$ | 2 | 3 | 4 | 5 | 6 | 7 | 8 | 9 | 10 | 11 | 12 | 13 | 14 | 15 | 16 |
|---|---|---|---|---|---|---|---|---|---|---|---|---|---|---|---|
| $f \leq$ | 35 | 21 | 21 | 23 | 26 | 28 | 32 | 35 | 37 | 40 | 43 | 46 | 49 | 51 | 55 |

Table 6: specified limits on $f$ for various $b$.

Our estimated Spaghetti boundary is sharp for sizes $n = 4\,(2),\,5\,(2),\,6\,(3),\,7\,(3)\,,8\,(4),\,9\,(5),\,12\,(7)$ (boundary between parentheses). For $n = 10, 11, 13, 15, 18, 21$ our boundary is sharp up to one unit. This involves some computer-provided facts related with question Q8 below. Note that improving the Spaghetti boundary may yield additional support for an earlier observation on thm. 4.1, that *most n-sets in a formal n-square can be transformed into an H-graph with just an H-move.*

**Q7.** *Find sharp Spaghetti boundaries for other sizes n.*

The problem of rotating a collection of 2-sets apart in a regular $n$-gon shows up in our attempts to improve the Spaghetti boundary of formal $n$-squares (see prop. 3.14). It can be shown that the regular $n$-gon (with $n$ even) admits a partition of its vertices into 2-sets assuming every possible diameter value if and only if $n$ is of type $8v$ or $8v+2$ (with $v>0$ integer). Allowing pairs with the same diameter complicates things. In a regular $n$-gon with $n = km$ and $k > 1$ odd, the pairs of diameter $m$ can be grouped into $m$ cycles of length $k$, whence no more than $\lfloor k/2 \rfloor \cdot m$ pairs can be disjoint.

Given any $n \geq 5$, define $b$ as $n/3$ if $n$ is divisible by 3 and $\lceil (n+2)/3 \rceil$ otherwise. For $n \leq 8$ we use lower rather than upper integer approximation. We verified with computer assistance that for $n \leq 38$ any set of $b$ pairs can be rotated apart. One easily verifies that this $b$ provides a *sharp* bound for $n \leq 11$ as well as for all $n$ divisible by 3 and for $n = 14, 20$ (cases $k = 3, 7, 5$ above). Exhaustive computer search revealed that for $n \leq 32$, these are the *only* sharp cases. Hence:

**Q8.** *What is the largest number b such that any b pairs in a regular n-gon can be rotated apart?*

The following observation on latin squares (Dénes end Keedwell [3, ch. 13, p. 467]) is a painful truth extending to equi-$n$-squares and perhaps to all problems raised above: *(...) for most computational problems, values of n greater than* 20 *are unmanageable.* The exact boundary may vary with the problem at hand. Our results being a mixture of reasoning and computing, to get around this obstacle requires increasing the share of reasoning --if possible.

# References

**[1]** R.A. Brualdi and H.J. Ryser, *Combinatorial Matrix Theory,* Cambridge University Press, Cambridge, UK, 1991.

**[2]** D. Bryant, J. Egan, B. Maenhaut and I.M. Wanless, *Indivisible Plexes in Latin Squares,* Des. Codes Cryptogr. 2009, DOI 10.1007/s10623-009-9269-z.

**[3]** J.H. Dénes and A.D. Keedwell, *Latin Squares and their Applications,* Academic Press, New York, 1974.

**[4]** J.H. Dénes and A.D. Keedwell (eds.), *Latin squares: New Developments in the Theory and Applications.* Annals of Discrete Mathematics, 46, (1991). Amsterdam: Academic Press. pp. xiv+454.

**[5]** J.H. Dénes and A.D. Keedwell, *Some Applications of Non-associative Algebraic Systems in*

*Cryptology,* P.U.M.A., Pure Math. Appl. 12 (2), (2001), 147-195.

**[6]** P. l'Ecuyer and R. Simard, *TestU01: A Software Library in ANSI C for Empirical Testing of Random Number Generators,* ACM Trans. Math. Software 33(4), Art. 22, (2007).

**[7]** J.B. Fraleigh, *A First Course in Abstract Algebra,* Pearson Education and Addison Wesley, Harlow, UK 2003, xii+520pp.

**[8]** O. Grosek and M. Sys, *Isotopy of Latin Squares in Cryptography,* Tatra Mt. Math. Publ., 45, (2010), 27-36.

**[9]** M.T. Jacobson and P. Matthews, *Generating Uniformly Distributed Random latin Squares,* J. Combinatorial Design 4 (1996), 405-437.

**[10]** David Joyner, *Adventures in Group Theory: Rubik's Cube, Merlin's Machine, and Other Mathematical Toys.* Johns Hopkins University Press, 2002. xviii + 262 pp.

**[11]** D.E. Knuth, *The Art of Computer Programming, Volume 2: Seminumerical Algorithms.* Addison-Wesley, Reading, Mass., 3rd edition, 1998. 762 pp.

**[12]** C. Koscielny, *NLPN Sequences over GF(q),* Quasigroups and Related Systems Vol. 4, 1997, pp. 89-102.

**[13]** R. Lidl, H. Niederreiter, *Introduction to Finite Fields and their Applications,* Cambridge University Press, 1994, 416pp.

**[14]** J.H. van Lint, R.M. Wilson, *A Course in Combinatorics.* Cambridge University Press, 2001 (2nd edition), 620pp.

**[15]** C.L. Liu, *Topics in Combinatorial Mathematics,* Mathematical Association of America, Buffalo, N.Y, 1972 (265pp).

**[16]** B.D. McKay and E. Rogoyski, *Latin Squares of Order Ten.* Electron. J. Combin., 2(3) (1995), 4 pp.

**[17]** B.D. McKay and I.A. Wanless, *On the Number of Latin Squares,* Annals of Combinatorics Vol. 9(3), (2005), pp 335-344.

**[18]** A.J. Menezes, P.C. van Oorschot, S.A. VanStone, *Handbook of Applied Cryptography,* CRC Press, Fifth Printing, 2001 (update 2010).

**[19]** K.A. Meyer, *A new Message Authentication Code Based on the Non-Associativity of Quasigroups,* Dissertation, Iowa State University, Ames, Iowa 2006, 91 pp.

**[20]** V.A. Shcherbacov *Quasigroups in Cryptology,* Computer Science Journal Moldova 17 No 2 (50), 2010, 194-217.

**[21]** S.K. Stein, *Transversals of Latin Squares and their Generalizations,* Pac. J. Math. 59(2), (1975), 567-575.

**[22]** M. Van de Vel, *Cryptologic Features of Shuffled Equi-n-Squares,* 2015, Preprint.

**[23]** D. Zwillinger (ed.), *Standard Mathematical Tables and Formulae,* (32nd edition) CRC Press, Boca Raton, Florida, 2011, 833pp.

## Appendix: The algorithm `Minrows`

A reader with Mathematica at hand may check some computational results in section 3 by entering the following code in a Mathematica notebook[6] (.nb file). The `Minrows` algorithm here is the "variant" one (see 3.5), returning the list of critical column partitions if the number of rows is critical (it returns the empty list or an excessive list otherwise). Note our special treatment of singletons in a partition (which is easily seen not to affect the result). Without this, Mathematica gets into problems with the dimension *n* approaching 40 (now postponed to $n \approx 50$).

```
rowunique[n_Integer, s_Integer, f_Integer] :=
 Max[Ceiling[s f/(n - s)], s - f] /; s < n

Minrows[n_Integer, r_Integer] := Module[
  {Crit = {}, K = {}, P = {}, PP = {}, plen, klen, rval, c1, cu, ru},
  For[m = r - 1, m > 1, m--,
   For[c = Ceiling[(n - m)/m], c <= r - m, c++,
    K = Select[IntegerPartitions[n - m, {c}], Max[#] <= m &];
    klen = Length[K];
    For[k = 1, k <= klen, k++,
     ru = Total[Map[rowunique[n, #, n - r] &, K[[k]]]];
     If[ru > r - m, Continue[]];
     P = Select[Join[{m}, K[[k]]], # > 1 &];(* drop singletons *)
     c1 = Length[P];
     cu = c + 1 - c1;(* number of singletons *)
     PP = Permutations[P];
     plen = Length[PP];
     For[p = 1, p <= plen, p++,
      rval = 0; (* check rows value *)
      For[j = 1, j <= c1, j++, rval += Ceiling[(n - rval) PP[[p, j]]/n]];
      If[rval + cu > r, Break[]]
      ];(* for p ... *)
     If[rval + cu <= r, AppendTo[Crit, Join[P, Table[1, {j, 1, cu}]]]]
     ] (* for k ... *)
    ](* for c ... *)
   ];(* for m ... *)
  Crit
  ](* end of module *)
```

---
**6**.        I used this code with Mathematica versions 7 and 8.